\documentclass[a4paper,11pt,reqno]{article}

\usepackage{amsmath,amsfonts,amssymb}
\usepackage[width=17cm,height=24cm]{geometry} 
\DeclareMathOperator{\expect}{{\mathbb E}}
\newcommand{\eps}{\varepsilon}
\newcommand{\ott}{[0,T]}
\makeatletter
\newcommand{\eqcolon}{\mathrel{\mathord{=}\raise.2\p@\hbox{:}}}
\newcommand{\coloneq}{\mathrel{\raise.2\p@\hbox{:}\mathord{=}}}
\makeatother
\newcommand{\der}{\delta}

\newcommand{\RR}{\mathbb{R}}
\newcommand{\TT}{\mathbb{T}}
\newcommand{\bbC}{\mathbb{C}}
\newcommand{\CC}{\mathcal{C}}

\newcommand{\DD}{\mathcal{D}}

\newcommand{\FF}{\mathcal{F}}
\newcommand{\QQ}{\mathcal{Q}}
\newcommand{\LL}{\mathcal{L}}
\newcommand{\YY}{\mathcal{Y}}

\newcommand{\ZZ}{\mathbb{Z}}

\newcommand{\norm}[1]{\lVert #1\rVert}

\newcommand{\wh}{\widehat}

\newcommand{\cZ}{\mathcal{Z}}
\DeclareMathOperator{\id}{\text{Id}}

\newtheorem{theorem}{Theorem}

\newtheorem{corollary}[theorem]{Corollary}

\newtheorem{lemma}[theorem]{Lemma}

\newtheorem{proposition}[theorem]{Proposition}
\newtheorem{remark}[theorem]{Remark}

\newenvironment{proof}[1][Proof]{\textbf{#1.} }{\ \rule{0.5em}{0.5em}}

\begin{document}
\title{Rough solutions for the\\ periodic Korteweg--de~Vries equation}
\author{
  { Massimiliano Gubinelli}              \\
{\small\it CEREMADE \& CNRS (UMR 7534) }\\[-0.1cm]
  {\small\it  Universit\'e Paris Dauphine}          \\[-0.1cm]
{\small \it  France}\\[-0.1cm] 
 {\small  {\tt gubinelli@ceremade.dauphine.fr}}   
}
\maketitle


\begin{abstract}
We show how to apply ideas from the theory of rough paths to the analysis of low-regularity solutions to non-linear dispersive equations. Our basic example will be the one dimensional Korteweg--de~Vries (KdV) equation on a periodic
domain and with initial condition in $\FF L^{\alpha,p}$ spaces.  We discuss
convergence of Galerkin approximations, a modified Euler scheme and the presence of a random force of white-noise type in time. \\[0.5cm]
\textbf{Keywords:} dispersive equations; rough paths; power series solutions.

\end{abstract}

In this note we start by considering the Cauchy problem for the
 classical   Korteweg--de~Vries (KdV) equation:
\begin{equation}
  \label{eq:kdv-real}
\partial_t u(t,\xi) + \partial^3_\xi u(t,\xi) +
\frac12 \partial_\xi u(t,\xi)^2 = 0,\quad u(0,\xi) = u_0(\xi), \qquad
(t,\xi)\in\RR\times \TT  
\end{equation}
where  the initial condition $u_0$ belongs to some Sobolev space $H^\alpha(\TT)$
of the torus $\TT = [-\pi,\pi]$. In a remarkable series of papers by Bourgain~\cite{MR1215780,MR1423623},
Kenig--Ponce--Vega~\cite{MR1329387,MR1230283,MR1086966} and later Colliander--Keel--Staffilani--Takaoka--Tao~\cite{MR1969209} this equation has been proved to possess global
solutions starting from initial conditions in $H^\alpha$ for any
$\alpha \ge -1/2$.   The existence of solutions in negative Sobolev spaces is possible due to the
regularizing effect of the dispersive linear term. This regularization
is more effective in the whole line and there global solutions exists for any $\alpha \ge -3/4$~\cite{MR1230283,MR1969209,Guo2009583}. Other references on the analysis of the KdV equation are~\cite{MR1019307,MR759907}.
More recently  Kappeler--Topalov, taking advantage of the complete integrability of this model, extended the global well-posedness in the periodic setting to any $\alpha \ge -1$ using the inverse scattering method~\cite{MR2181576}.  

Inspired by the theory of \emph{rough paths} we will look for an alternative approach to the
construction of solutions of eq.~(\ref{eq:kdv-real}) and in general for dispersive equations with polynomial non-linearities. Our method turns out to have some similarities with Christ's power series approach~\cite{MR2333210, MR2411067} and to allow to consider KdV with initial condition in $\FF L^{\alpha,p}$ for $1\le p \le +\infty$.

Rough paths (introduced by Lyons in~\cite{MR1654527,MR2036784}) allow the study of differential
equations driven by irregular functions. They have been applied to the
path-wise study of stochastic differential equations driven by
Brownian motion, by fractional Brownian motion of any index $H>1/4$~\cite{MR1883719}
and other stochastic processes~\cite{MR2164030,fvsubell}. Part of the theory has been reformulated in terms of the sub-Riemmannian geometry of certain Carnot groups~\cite{fvgeo}. 
In~\cite{MR2091358} we showed how to reinterpret the work of Lyons in the terms of a
cochain complex of finite increments and a related integration
theory. The key step is the introduction of a map $\Lambda$ (called
sewing map) which encodes a basic fact of rough path theory.
We exploited this point of view 
 to treat stochastic partial differential
equations of evolution type~\cite{youngspde,TindelGubinelli} and to
study the initial value problem for a partial differential  equation
modeling the approximate evolution of random vortex filaments in 3d
fluids~\cite{evolution}. 

We would like to show that the concepts of the theory can be used
fruitfully for problems not related to stochastic processes.
The periodic Korteweg--de Vries equation is used as a case study for our
ideas. The point of view here developed should be applicable also to KdV on
the full  real line or other dispersive semi-linear equations like the
modified KdV or the non-linear Schr\"odinger equation (indeed the results~\cite{MR2333210, MR2411067} obtained via the power-series method can be understood in terms of rough paths following the lines of the present investigation).  

No previous knowledge of rough paths theory is required nor any
 result on the periodic KdV is used in the following. We took care to
make the paper, as much as possible, self-contained. Since the
arguments are similar to those used in finite dimensions
in~\cite{MR2091358},  the reader can refer to this last paper to gain a wider perspective
on the technique and on the stochastic applications of rough paths. 

Let us describe the main results of this paper in terms of distributional solutions of KdV:
let $P_N$ be the Fourier projector on modes $k\in \ZZ$ such that $|k|\le N$ and $\mathcal{N}(\varphi)(t,\xi)=\partial_\xi (\varphi(t,\xi)^2)/2$ for smooth functions $\varphi$. 
 Then
\begin{theorem}
\label{th:main-intro}
For any $1 \le p \le \infty$, $\alpha > \alpha_{*}(p)=\max(-1/p,-1/2)$ and $u_0 \in \FF L^{\alpha,p}$ there exists a $T_* >0$ and a continuous function $u(t) \in C([0,T^*],\FF L^{\alpha,p})$ with $u(0) = u_0$ for which the distribution $\mathcal{N}(P_N u)$ converge as $N\to \infty$ to a limit which we denote by 
$\mathcal{N}(u)$ and moreover the distributional equation
$$
\partial_t u + \partial_\xi^3 u +  \mathcal{N}(u) = 0
$$
is satisfied in $[0, T^*]\times \TT$. 
\end{theorem}
This solution is the limit of smooth solutions and of some modified  Galerkin approximations. There exists a natural space of continuous functions on $\FF L^{\alpha,p}$ for which the nonlinear term can be defined as a distribution and where uniqueness of solutions holds. 
We show also how to implement the $L^2(\TT)$ conservation law in the rough path approach and obtain in this way global solutions in $L^2(\TT)$. 
Following the lines of the numerical study of SDEs driven by rough paths~\cite{davie,MR1470933,fveuler} we analyze an Euler-like time-discretization of the PDE which converges to the above solution. 
Finally we also prove an existence and uniqueness result under random perturbation of white-noise type in time of the form
\begin{equation*}
\partial_t u + \partial^3_\xi u +
\frac12 \partial_\xi u^2 = \Phi \partial_t  \partial_\xi B
\end{equation*}
where $\Phi$ is a bounded linear operator from $\FF L^{0,\infty}$ to $\FF L^{\alpha,p}$ and
$\partial_t \partial_\xi B$ a white noise on $\RR\times\TT$.

\bigskip
\paragraph{Plan.}
 In Sect.~\ref{sec:intro} we start by recasting the KdV
eq.~(\ref{eq:kdv-real}) in its mild form and to perform some
manipulation to motivate the finite-increment equation which we will study
in Sect.~\ref{sec:rough}  where we prove  existence and uniqueness of
local solutions and discuss the distributional meaning of these solutions.
 Then we prove a-priori estimates necessary
show that global rough solutions
exists for initial conditions in $L^2(\TT)$ (Sect.~\ref{sec:global}).
In
Sect.~\ref{sec:galerkin} we  prove that the rough solutions are
limits of suitable Galerkin approximations. Moreover in Sect.~\ref{sec:euler} we introduce a
time discretization scheme and prove its convergence.
 The
equation that we study is just one of a stack of finite-increment equations that can be
generated starting from KdV. In Sect.~\ref{sec:more} we derive the
next member of this hierarchy. Finally
Sect.~\ref{sec:stoch} addresses the problem of the presence of an
additive stochastic forcing.
We collect some longer and technical proofs in the Appendix.

\paragraph{Notations.}
 We denote with
$\wh f: \ZZ \to \bbC$ the Fourier coefficients of a real function $f : \TT \to
\RR$:
$
f(\xi) = \sum_{k \in \ZZ} \hat f(k) e^{i k \xi}
$
and define the space $$
\FF L^{\alpha,p} = \{ f \in \mathcal{S}'(\TT) : \hat f(0)=0, |f|_{\FF L^{\alpha,p}}=|\langle \cdot \rangle^\alpha \hat f |_{\ell^p} < \infty \}
$$
where $ \mathcal{S}'(\TT)$ is the space of real Schwartz distributions on the torus $\TT$ and $\langle x \rangle = (1+|x|^2)^{1/2}$. We restrict the space to the mean zero functions since this is the natural setting to discuss periodic KdV.
 Then $H^\alpha(\TT)\backslash \RR = \FF L^{\alpha,2}$.  Sometimes we note $\FF$ the Fourier transform operator so that $\FF f  = \hat f$. Given two Banach spaces $V,W$, denote with $\LL(V,W)$ the Banach space of bounded linear operators from $V$ to $W$ endowed with the operator norm and with $\LL^n V = \LL(V^n, V)$. If $X \in \LL^n V$ we  write $X(a) = X(a,\cdots,a)$ when the operator is applied to $n$ copies of the same argument $a\in V$.  
The symbol $C$ in the r.h.s. of  estimates denotes a positive constant which can be different from line to line and usually we write $A \lesssim B$ for $A \le C B$.

\bigskip

\section{Formulation of the problem}
\label{sec:intro}
Duhamel's perturbation formula applied to eq.~(\ref{eq:kdv-real}) gives
\begin{equation}
\label{eq:kdv-mild}
u(t) = U(t) u_0 + \frac{1}{2} \int_0^t U(t-s) \mathcal{N}(u(s))\, ds 
\end{equation}
where  $U(t)$ is the Airy group of isometries of any $\FF L^{\alpha,p}$ given by the
solution of the linear part of eq.~(\ref{eq:kdv-real}). The group $U(t)$ acts on $\varphi$ as
$
\FF (U(t) \varphi)(k) = e^{i k^3 t} \wh \varphi(k)
$.
 Conservation of mass for KdV guarantees that  if
$\wh u_0(0) = 0$ then we will have $\wh u(t,0) = 0$ for any $t \ge 0$.
Using as
unknown the twisted variable $v(t) = U(-t) u(t)$ we have
\begin{equation}
  \label{eq:kdv-base}
 v(t) =   v_0 + \int_0^t \dot X_s( v(s))\, ds , \quad t \in [0,T]\end{equation}
where $v_0 = u_0$, $\dot X_s( v(s)) = \dot X_s( v(s),v(s))$ and where the bilinear unbounded operator $\dot X_s \in \LL^2 \FF L^{s,p}$ is defined as $\dot X_s (\varphi_1, \varphi_2) = U(-s)\partial_\xi[(U(s)\varphi_1) (U(s)\varphi_2)]/2$ so that its Fourier transform reads
$$
\FF \dot X_s(\varphi_1,\varphi_2)(k) = 
\frac{ik}2 \sum'_{k_1} e^{- i 3 k k_1 k_2 s}
  \wh \varphi_1(k_1) \wh \varphi_2(k_2) 
$$
with the convention that $k_2 = k-k_1$ and where we use the
convention that the primed summation over $k_1 \in \ZZ$ does not include the
terms when $k_1 = 0$ or $k_1 = k$. In obtaining this equation
we have used the algebraic identity $k^3 - k_1^3 - k_2^3 = 3 k k_1
k_2$ valid for any triple $k,k_1,k_2$ such that $k = k_1 + k_2$.
In the rest of the paper we will concentrate on the study of solution
to eq.~(\ref{eq:kdv-base}) with the fixed point method. 

\subsection{Power series solutions and generalized integration}

Proceeding formally we can expand any solution to~\eqref{eq:kdv-real} into a series involving only the operators $\dot X$ and the initial condition:
\begin{equation}
\label{eq:power-series}
\begin{split}
v(t) = v_0 + \int_{0\le \sigma\le t} \dot X_{\sigma}(v_0,v_0) d\sigma+ 2 \int_{0\le \sigma_1\le\sigma\le t}  \dot X_{\sigma}(v_0,\dot X_{\sigma_1}(v_0,v_0)) d\sigma d\sigma_1
\\
+  4 \int_{0\le \sigma_2 \le \sigma_1\le\sigma\le t}  \dot X_{\sigma}(v_0,\dot X_{\sigma_1}(v_0,\dot X_{\sigma_2}(v_0,v_0))) d\sigma d\sigma_1 d\sigma_2
\\
+  2 \int_{0\le \sigma_{1,2} \le \sigma\le t}   \dot X_{\sigma}(\dot X_{\sigma_2}(v_0,v_0),\dot X_{\sigma_1}(v_0,v_0)) d\sigma d\sigma_1 d\sigma_2
+ \cdots 
\end{split}
\end{equation}
This perturbative expansion is naturally indexed by (binary) trees representing the various ways of applying $\dot X$ to itself and of performing time integrations. A possible approach to define a solution to the differential equation is then to prove that each term of the series is well defined and that the sum of the series converges. There is quite a bit of literature on this method for Navier-Stokes like equations~\cite{Gall,LS,Sinai3,Sinai2,MR2227041} and recently Christ~\cite{MR2333210} advocated this approach in the context of non-linear dispersive equations and used it to solve a modified periodic non-linear Schr\"odinger equation in  $\FF L^{\alpha,p}$ spaces. See also~\cite{MR2411067} for another recent paper studying power-series solution for the periodic modified KdV equation using Christ's approach.
 
The advantage of the power series expansion 
is that the relevant objects are the operators obtained from $\dot X$ by successive application and time integrations and applied to the inital condition. For example the first term in the expansion involves the the symmetric bilinear operator
$X_{ts}$ defined  for any $s\le t$ as
$$
\FF X_{ts}(\varphi_1,\varphi_2)(k) = \int_s^t \FF \dot X_\sigma(\varphi_1,\varphi_2)(k) d\sigma = 
\sum'_{k_1} \frac{e^{-i 3 k k_1 k_2
      s}-e^{-i 3 k k_1 k_2 t}}{6 k_1 k_2}\wh\varphi_1(k_1)
\wh \varphi_2(k_2) .
$$
While $\dot X_\sigma$ is usually unbounded in any $\FF L^{\alpha,p}$, the operator $X_{ts}$ will be shown to be bounded in $\FF L^{\alpha,p}$ for appropriate values of $\alpha,p$ and a similar  behavior will hold for higher order operators. The regularizing effect of time integrations is due to the dispersive nature of the linear part combined with absence of resonances in the non-linear term.
  
  On the other end the main disadvantage  of the power-series method is that one has to control  arbitrary terms of the series and the proof of summability require usually a big analytical and combinatorial effort. 
  Another interesting difference is that in the rough path approach there is available a natural space where uniqueness of solutions holds while in the power-series approach no uniqueness result is available.

  Rough path theory allows to bypass the complete power-series expansion exploiting the smallness  of the leading terms for small time intervals and using a generalized notion of integration to pass from the approximate solution for an very small time interval to an exact solution for an $O(1)$ interval of time. 
Instead of writing the series solution from time $0$ to time $t$ we write it between two times $ s \le t$ denoting by $\der v_{ts} = v_t - v_s$ with $v_t = v(t)$ (the index
notation being more comfortable in the following) so that
\begin{equation}
\label{eq:kdv-base-remainder}
 \der v_{ts} = X_{ts}(v_s) + r_{ts}
\end{equation}
where $r_{ts}$ stands for the rest of the series and will be treated  as a negligible remainder term. For this to make sense we need that $X_{ts}$ gives itself a small contribution when $|t-s|$ is small. Our first result give then a quantitative control of the size of $X_{ts}$ as a bounded operator in $\FF L^{\alpha,p}$.
When $2\le p\le +\infty$ define the set $\DD\subset \RR \times \RR_+$ of pairs $(\gamma,\alpha)$ by
$$\DD =\{
\alpha \ge -\frac{1}{2}-\frac 1 p + \gamma, 0 \le \gamma < \frac14 
\} \cup
\{
\alpha > -1-\frac 1 p + 3\gamma, \frac1 4 \le \gamma \le \frac12 \} ;
$$
while when $1 \le p \le 2$ 
$$
\DD = \{
\alpha \ge -1 + \gamma, 0 \le \gamma < \frac1{2p} \} \cup
\{
\alpha > -1-\frac 1 p + 3\gamma, \frac1 {2p} \le \gamma \le \frac12\} .
$$
Then
\begin{lemma}
\label{lemma:X}
For any couple
$(\gamma,\alpha) \in \DD$ 
the operator $X_{ts}$ is bounded from $(\FF L^{\alpha,p})^2$ to
$\FF L^{\alpha,p}$ and
$
|X_{ts}|_{\LL^2 \FF L^{\alpha,p}} \lesssim_{\gamma,\alpha} |t-s|^\gamma .
$
\end{lemma}

The proof of this lemma is postponed to Appendix~\ref{sec:reg}.
  A trivial but important observation is that the family of operators $\{X_{ts} \}_{ s \le t }$ satisfy the algebraic equation
\begin{equation}
  \label{eq:X-chain}
X_{ts}(\varphi_1,\varphi_2) - X_{tu}(\varphi_1,\varphi_2) -
X_{us}(\varphi_1,\varphi_2) = 0, \quad  s\le u\le t 
\end{equation}
as can be easily checked using the definition.

\subsection{The sewing map}

The first term in the r.h.s. of eq.~(\ref{eq:kdv-base-remainder})  is well understood thanks to
Lemma~\ref{lemma:X}, while the $r$ term contains all the difficulty. However, due to the particular
structure of eq.~(\ref{eq:kdv-base-remainder}), the $r$ term must satisfy a simple
algebraic equation. Indeed, for any triple
$s\le u\le t$ we have
$
\der v_{ts} - \der v_{tu} - \der v_{us} = 0
$ and substituting eq.~(\ref{eq:kdv-base-remainder}) in this relation we get
\begin{equation}
\label{eq:pre-lambda}
  \begin{split}
 r_{ts} -
r_{tu} - r_{us}   = 
- X_{ts}(v_s,v_s) + X_{tu}(v_u,v_u) + X_{us}(v_s,v_s).   
    = 
 X_{tu}(\der v_{us},v_s) + X_{tu}(v_u,\der v_{us})   
  \end{split}
\end{equation}
where we used eq.~(\ref{eq:X-chain}) to simplify the r.h.s..
The main observation contained in the work~\cite{MR2091358} is that  sometimes this
equation determines $r$ uniquely.
To explain the conditions under which we can solve
eq.~(\ref{eq:pre-lambda}) we need some more notation. 
Given a normed vector space $(V,|\cdot|)$ introduce the vector space
$\CC_n V \subset C([0,T]^n;V)$ such that $a \in \CC_n V$ iff $ a_{t_1, \dots,t_n} = 0$ when $t_i
  = t_j$ for some $1 \le i < j \le n$. We have already introduced the
  operator $\der : \CC_1 V \to \CC_2 V$ defined as $\der f_{ts} = f_t
  - f_s$. Moreover, is useful to introduce
another operator $\der : \CC_2 V \to \CC_3 V$ defined on
continuous functions of \emph{two} parameters 
on a vector space $V$ as
$ \der a_{tus} = a_{ts} - a_{tu} - a_{us} $. The two operators satisfy
the relation $\der \der f = 0$ for any $f \in \CC_1 V$. Moreover if $a
\in \CC_2 V$ is such that $\der a = 0$ then  there exists $f \in \CC_1 V$ such that $ a = \der f$. Denote
$\cZ \CC_3 V = \CC_3 V \cap \mathrm{Im} \der$, where $\mathrm{Im} \der$ is the image of the $\der$ operator.
We measure the size of elements in $\CC_n V$ for $n=2,3$  by H\"older-like norms
defined in the following way: for $f \in \CC_2 V $ let
$$
\norm{f}_{\mu} =
\sup_{s,t\in\ott}\frac{|f_{ts}|}{|t-s|^\mu},
\quad\mbox{and}\quad
\CC_1^\mu V =\{ f \in \CC_2 V;\, \norm{f}_{\mu}<\infty  \}.
$$
In the same way, for $h \in \CC_3 V$, set
\begin{eqnarray}
  \label{eq:normOCC2}
  \norm{h}_{\gamma,\rho} &=& \sup_{s,u,t\in\ott} 
\frac{|h_{tus}|}{|u-s|^\gamma |t-u|^\rho}\\
\|h\|_\mu &= &
\inf\left \{\sum_i \|h_i\|_{\rho_i,\mu-\rho_i} ;\, h =
 \sum_i h_i,\, 0 < \rho_i < \mu \right\} ,\nonumber
\end{eqnarray}
where the last infimum is taken over all finite sequences $\{h_i \in \CC_3 V  \}$ such that $h
= \sum_i h_i$ and for all choices of the numbers $\rho_i \in (0,z)$.
Then  $\|\cdot\|_\mu$ is easily seen to be a norm on $\CC_3 V$, and we set
$
\CC_3^\mu V =\{ h\in\CC_3 V ;\, \|h\|_\mu<\infty \}
$.
Eventually,
let $\CC_n^{1+} V  = \cup_{\mu > 1} \CC_n^\mu V $ for $n =2,3$
and remark that the same kind of norms can be considered on the
spaces $\cZ \CC_3 V $, leading to the definition of the spaces
$\cZ \CC_3^\mu V $ and $\cZ \CC_3^{1+} V $. 
The following proposition is the basic result which allows the
solution of equations in the form~(\ref{eq:pre-lambda}).
\begin{proposition}[The sewing map $\Lambda$]
\label{prop:Lambda}
There exists a unique linear map $\Lambda: \cZ \CC^{1+}_3 V 
\to \CC_2^{1+} V $ such that 
$
\delta \Lambda  = \id_{\cZ \CC_3 V }
$.
Furthermore, for any $\mu > 1$, 
this map is continuous from $\cZ \CC^{\mu}_3 V $
to $\CC_2^{\mu} V $ and we have
\begin{equation}\label{ineqla}
\|\Lambda h\|_{\mu} \le \frac{1}{2^\mu-2} \|h\|_{\mu} ,\qquad h \in
\cZ \CC^{1+}_3 V . 
\end{equation}
\end{proposition}
This proposition has been first proved in~\cite{MR2091358}. A simplified proof is contained in~\cite{TindelGubinelli}.
Using the notations just introduced, eq.~(\ref{eq:pre-lambda}) take
the form
\begin{equation}
  \label{eq:r-1}
\der r_{tus} = X_{tu}(\der v_{us},v_s) + X_{tu}(v_u, \der v_{us}) = [X(\der v,v) + X(v, \der v)]_{tus},
\end{equation}
where, by construction, the r.h.s. belongs to $\cZ \CC_3 \FF L^{\alpha,p}$. If
we can prove that  it actually belongs to $ \cZ \CC^z_3 \FF L^{\alpha,p}$
for some $z >1$, Prop.~\ref{prop:Lambda} will give us the possibility to state that the
unique solution of eq.~(\ref{eq:r-1}) in $\CC^{1+}_2 \FF L^{\alpha,p}$ is given by
$r = \Lambda [X(\der v,v) + X(v, \der v)]$. 

\subsection{$\Lambda$ equations}

Let us come back to our initial problem. Since we aim to work in distributional spaces the rigorous meaning of eq.~(\ref{eq:kdv-base}) is a priori not clear (even in a weak sense). By formal manipulations we have been able to recast the initial problem in a finite-difference equation involving the $\Lambda$ map which reads:
\begin{equation}
  \label{eq:abs-lambda}
\der v = X(v,v) + \Lambda[X(\der v, v) + X(v, \der v)]   
\end{equation}
where we used an abbreviated notation since all the terms have been
already described in detail. Note that in this equation the $\Lambda$ map has replaced the integral
 so we
will call this kind of equations: \emph{$\Lambda$-equations}. Instead
of solving the integral equation 
we would like to  solve, by fixed-point
methods, the $\Lambda$-equation~(\ref{eq:abs-lambda}). Afterwards we will show rigorously that solution to such $\Lambda$-equations gives generalized solutions to KdV.

  Unfortunately,
for the particular form of our $X$ operator, we are not able to show
that this equation has solutions. Recall  Lemma~\ref{lemma:X} where we showed that $X$ belongs to 
$\CC_2^\gamma \LL^2 \FF L^{\alpha,p}$ for $(\gamma,\alpha) \in \DD$ so we
must expect that $\der v \in \CC_2^\gamma  \FF L^{\alpha,p}$ since the
$\Lambda$ term will belong (at worst) to $\CC_2^{1+}  \FF L^{\alpha,p}$. But then we have
$
X(\der v,v) + X(v,\der v) \in \CC_3^{2\gamma} \FF L^{\alpha,p}
$ 
so that we must require $2\gamma > 1$ in order for this term to be in
the domain of $\Lambda$. The  set $\DD$ however contains only
values of $\gamma \le 1/2$ so we are not able to study
eq.~(\ref{eq:abs-lambda}) in $\CC_1^\gamma  \FF L^{\alpha,p}$ for
$(\gamma,\alpha) \in \DD$.

This difficulty can be overcome by truncating the power-series solution~\eqref{eq:power-series}
to higher order. By looking at
eq.~(\ref{eq:kdv-base-remainder}) we see that the next step in the
expansion will generate the following expression
\begin{equation}
  \label{eq:exp-step-2}
  \begin{split}
 \der v_{ts} & =   X_{ts}(v_s)+X^2_{ts}(v_s)  + r^{(2)}_{ts}
  \end{split}
\end{equation}
where the  trilinear operator  $X^2 \in \CC_2 \LL^3 \FF L^{\alpha,p}$ is defined by
\begin{equation}
  \label{eq:X2}
  \begin{split}
 X^2_{ts}&(\varphi_1,\varphi_2,\varphi_3) = 2 \int_s^t d\sigma \int_s^{\sigma_1}d\sigma_1 \, \dot X_\sigma(\varphi_1,  \dot X_{\sigma_1}(\varphi_2, \varphi_3))
 \end{split}
\end{equation}
and satisfies the algebraic relation
\begin{equation}
  \label{eq:area-eq}
  \der X^2(\varphi_1,\varphi_2,\varphi_3)_{tus} = 2 X_{tu}(\varphi_1,X_{us}(\varphi_2,\varphi_3))
\qquad s\le u \le t
\end{equation}
which show that $X^2$ is related to the previously defined operator
$X$ by a ``quadratic'' relationship (this motivates the abuse of the superscript $2$
in the
definition of $X^2$).  
For the regularity of $X^2$ we have the following result whose proof  is again postponed to Appendix~\ref{sec:reg}.
\begin{lemma}
\label{lemma:X2}
There exists unbounded operators $\hat X^2, \breve X^2 \in \LL^3\FF L^{\alpha,p}$  such that $X^2_{ts}= \hat X^2_{ts}+\breve X^2_{ts}$, $\der \breve X^2 = 0$ and when $\alpha > \alpha_{*}(p)=\max( -1/p,-1/2)$ we have $|\breve X_{ts}|_{\LL^3 \FF L^{\alpha,p}} \lesssim_{\alpha,\gamma} |t-s|$ and for any couple
$(\gamma,\alpha) \in  \DD$
we have
$
|\hat X^2_{ts}|_{\LL^3 \FF L^{\alpha,p}} \lesssim_{\alpha,\gamma} |t-s|^{2\gamma} 
$.
\end{lemma}
Note that the full operator $X^2$ can be controlled only in the smaller
region $\DD'= \DD \cap \{\alpha> \alpha_{*}(p)\}$ of the $(\gamma,\alpha)$ plane. This limitation has some
connection with the low regularity ill-posedness discussed in~\cite{MR2018661} and is related to the periodic setting. 
To infer the $\Lambda$-equation  related to eq.~(\ref{eq:exp-step-2}) we simply apply to it $\der$ operator and we get
\begin{equation}
  \label{eq:abs-step-2-ur-rem}
  \begin{split}
\der \hat r^{(2)}_{tus} = & 2 X_{tu}(\der v_{us}- X_{us}(v_s,v_s),v_s) 
+ X_{tu}(\der v_{us},\der v_{us}) 
\\ & + X^2_{tu}(\der v_{us},v_s,v_s) +
X^2_{tu}(v_u,\der v_{us},v_s) + X^2_{tu}( v_u,v_u,\der v_{us}).     
  \end{split}
\end{equation}

\section{Rough solutions of KdV}
\label{sec:rough}
From eq.~(\ref{eq:abs-step-2-ur-rem}) we can write down
 a second $\Lambda$-equation (beside eq.~(\ref{eq:abs-lambda})) associated to the KdV equation:
\begin{equation}
  \label{eq:abs-step-2}
\der v = (\id - \Lambda \der)[X(v)+ X^2(v)] = X(v) +  X^2(v) + v^\flat  
\end{equation}
with
$
v^\flat = \Lambda[ 2 X(\der v- X(v),v) 
+ X(\der v,\der v) 
 +  X^2(\der v,v,v) +
 X^2(v,\der v,v) +  X^2( v,v,\der v)   ]
$.
To give a well-defined meaning to eq.~(\ref{eq:abs-step-2}) it will be enough that the argument of $\Lambda$ belongs actually to its domain. Given an allowed pair $(\gamma,\alpha) \in \DD'$ which fixes the regularity of $X$ and $ X^2$ sufficient (and natural) requirements for $v$ are
\begin{equation}
  \label{eq:requirements}
\sup_{t\in[0,T]}|v_t|_{\FF L^{\alpha,p}} <\infty,\qquad  \der v \in \CC^\gamma_1 \FF L^{\alpha,p},\qquad \der v - X(v,v) \in \CC^{2\gamma}_2 \FF L^{\alpha,p}.   
\end{equation}
Under these conditions $X(\der v- X(v,v),v) \in \CC^{3\gamma}_2  \FF L^{\alpha,p}$ and also $ X^2(\der v,v,v) +
 X^2(v,\der v,v) +  X^2( v,v,\der v) \in \CC^{3\gamma}_2 L^{\alpha,p}$. They are in the domain of $\Lambda$ if $3\gamma > 1$ i.e. $\gamma > 1/3$. This fixes the limiting time regularity for the $\Lambda$ eq.~(\ref{eq:abs-step-2}) and it turns out that for any $1\le p \le \infty$ and $\alpha > \alpha_*(p)$ there is a pair $(\gamma,\alpha) \in \DD'$ with $\gamma > 1/3$. This will fix the regularity of the initial data that we are able to handle.
 
To define solutions of eq.~(\ref{eq:abs-step-2}) let us introduce a suitable  space to enforce the all the conditions in eq.~(\ref{eq:requirements}).
For any $0 < \eta \le \gamma$ and consider the complete metric space $\QQ_\eta$ whose elements are triples $(y,y',y^\sharp)$
where $y \in \CC_1^\eta  \FF L^{\alpha,p}$, $y' \in \CC_1^\eta  \FF L^{\alpha,p}$, $y'_0 = y_0$ and $y^\sharp \in \CC_2^{2\eta}
 \FF L^{\alpha,p}$. Additional requirement is that they satisfy the equation $\der y = X(y',y') +
y^\sharp$.  The distance $d_{\QQ,\eta}$ on $\QQ_\eta$ is defined by
$$
d_{\QQ,\eta}(y,z) = |y_0-z_0| + \|y-z\|_{\eta} + \|y'-z'\|_{\eta} + \|y^\sharp-z^{\sharp}\|_{\eta}. 
$$
for any two elements $y,z \in \QQ_\eta$. 
With abuse of notations we will denote a triple
$(y,y',y^\sharp) \in \QQ_\eta$ using simply its first component.
Moreover  we denote with $y_0$ also the constant path in $\QQ_\eta$ with value $(y_0,0,0)$. 
The main result of the paper is the following theorem.

\begin{theorem}
\label{thm:main}
For any $\alpha > \alpha_*(p)$ take $\gamma > 1/3$ such that $(\gamma,\alpha) \in \DD'$.
Then, for any
$v_0 \in  \FF L^{\alpha,p}$ and for a sufficiently small interval of time $[0,T^*]$ where $T^*$ depends only on the norm of $v_0$, there exists a unique $v \in \CC^\gamma \FF L^{\alpha,p}$ such that $v(0) = v_0$ and 
$$
v_t = v_s + X_{ts}(v_s) +  X^2_{ts}(v_s) + o(|t-s|)
$$
for all $0\le s\le t \le T^*$.
If we write $v = \Theta(X, X^2,v_0)$ then  
$
\Theta: \CC^\gamma_2 \LL^2 \FF L^{\alpha,p} \times \CC^{2\gamma}_2 \LL^3\FF L^{\alpha,p} \times  \FF L^{\alpha,p} \to \QQ_\gamma
$  is a locally Lipschitz map wrt all its arguments.
\end{theorem}
Since $(\gamma,\alpha)\in\DD$ the operators $X$ and $ X^2$ are regular enough so that the proof of the theorem can
follow essentially the pattern of a similar result
in~\cite{MR2091358}. A proof is given in Appendix~\ref{proof:thm:main}.
Since $v\in\QQ_\gamma$ we have that $v^\flat \in \CC_2^{3\gamma}\FF L^{\alpha,p}$ with $3\gamma > 1$, then $v$ satisfy the property
\begin{equation}
\label{eq:rough-sum}
v_t = v_0 + \lim_{|\Pi|\to 0} \sum_{i} X_{t_{i+1} t_i}(v_{t_i}) +  X^2_{t_{i+1} t_i}(v_{t_i})
\end{equation}
where $\Pi = \{0=t_0\le t_1 \le \cdots \le t_n =t \}$ is a partition of $[0,t]$ and $|\Pi| = \max |t_{i+1}-t_i|$ its size. The limit is in $\FF L^{\alpha,p}$. The proof is simple: by eq.~\eqref{eq:abs-step-2} 
$$
\sum_{i=0}^{n-1} X_{t_{i+1} t_i}(v_{t_i}) +  X^2_{t_{i+1} t_i}(v_{t_i}) = \sum_{i=0}^{n-1} (\der v)_{t_{i+1} t_i} - \sum_{i=0}^{n-1}  v^\flat_{t_{i+1} t_i} = v_t - v_0 - \sum_{i=0}^{n-1}  v^\flat_{t_{i+1} t_i} 
$$
and the second term goes to zero as $|\Pi| \to 0$. Below we will also prove that an Euler scheme related to eq.~\eqref{eq:rough-sum} converges as the size of the partition goes to zero. A nice property of these solutions is the following (already  noted by Christ for power series solutions of NLS~\cite{MR2333210}).

\begin{corollary}
Let $v$ be the unique solution of the $\Lambda$-equation given by Thm.~\ref{thm:main} and let $u(t) = U(t) v(t)$. Let $P_N$ be the Fourier projector on modes $k$ such that $|k|\le N$ and let $\mathcal{N}(\varphi)(t,\xi)=\partial_\xi (\varphi(t,\xi)^2)/2$ for smooth functions $\varphi$. Then  in the sense of distributions in $\CC([0,T^*], \mathcal{S}'(\TT))$ we have convergence of $\mathcal{N}(P_N u)$ to a limit which we denote by 
$\mathcal{N}(u)$ and moreover the distributional equation
$$
\partial_t u + \partial_\xi^3 u +  \mathcal{N}(u) = 0
$$
is satisfied.
\end{corollary}

\begin{proof}
We start by proving that $\mathcal{N}(P_N u) \to \mathcal{N}(u)$ distributionally. It is enough to prove for any $0\le s \le t \le T^*$ the convergence of 
$
V_{t}=\int_0^t  U(-r)\mathcal{N}(P_N u(r)) dr 
$
in $\FF L^{\alpha,p}$ since any smooth test function can be approximated in time by step functions. Now using eq.~\eqref{eq:abs-step-2} we have
\begin{equation}
\label{eq:weak-conv-decomp}
\der V_{ts}=
\int_s^t \dot X_r(P_N v(r),P_N v(r)) dr
=Y_{ts}(v_s,v_s)+Y^2_{ts}(v_s,v_s,v_s) + V^\flat_{ts}
\end{equation}
where $V^\flat$ is a remainder term, $Y_{ts} = X_{ts}(P_N \times P_N)$ and
$$
Y^2_{ts}(\varphi_1,\varphi_2,\varphi_3) = 2 \int_s^t \dot X_r(P_N \varphi_1,P_N X_{rs}(\varphi_2,\varphi_3)) dr .
$$
From the regularity proofs for $X$ and $X^2$ the following facts are easy to prove:
 (i)
 $Y, Y^2$ enjoy at least the same regularity of $X$ and $ X^2$; (ii)
as $N\to \infty$ they are equibounded in $\CC^\gamma \LL^2 \FF L^{\alpha,p}$ and $\CC^{2\gamma} \LL^3 \FF L^{\alpha,p}$ respectively; (iii) $Y_{ts}\to X_{ts}$ and $ Y^2_{ts}\to  X_{ts}^2$ in the strong operator norm for fixed $t,s$.
For $V^\flat_{ts}$ we have then the following equation
$$
\der V^\flat_{tus} = 2 Y_{tu}(\der v_{us}-X_{us}(v_s,v_s),v_s)+Y_{tu}(\der v_{us},\der v_{us})
$$
$$
+  Y^2_{tu}(\der v_{us},v_s,v_s)+ Y^2_{tu}(v_u,\der v_{us},v_s) +  Y^2(v_u,v_u,\der v_{us})
$$
(cfr. eq.~\eqref{eq:abs-step-2-ur-rem}).
  Using the sewing map we have that the functions $V^\flat$ are also equibounded in $\CC_3^{3\gamma}\FF L^{\alpha,p}$. For fixed $0\le s\le t \le T^*$ we have
$$
\der V_{ts}-V^\flat_{ts} = Y_{ts}(v_s)+ Y^2_{ts}(v_s) \to 
X_{ts}(v_s)+ X^2_{ts}(v_s) = \der v_{ts}-v^\flat_{ts}
$$
so $\limsup_{N\to \infty}\ |\der(V-v)_{ts} |_{\FF L^{\alpha,p}} = \limsup_{N\to \infty} | V^\flat_{ts}-v^\flat_{ts}|_{\FF L^{\alpha,p}}†\lesssim |t-s|^{3\gamma} 
$
and since $3\gamma > 1$ this implies that $\limsup_{N\to \infty} \sup_{0\le t\le T^*} |V_t-(v_t-v_0)|_{\FF L^{\alpha,p}} = 0$ proving that $V \to v-v_0$ in $\CC^0_1 \FF L^{\alpha,p}$. This gives the distributional convergence of $\mathcal{N}(P_N u)$. If we call $\mathcal{N}(u)$ the limit we have $V_t = \int_0^t U(-r)\mathcal{N}(u(r)) dr$ and
$
u(t) = U(t)u(0) + \int_0^t U(t-r)\mathcal{N}(u(r)) dr
$
which is the mild form of the differential equation.
\end{proof}

\medskip

Some remarks are in order. The function $v$ is $\gamma$-H\"older continous in $\FF L^{\alpha,p}$,  $|u(t)-u(s)|_{\FF L^{\alpha,p}} \le |(U(t)-U(s))u(s)|_{\FF L^{\alpha,p}} + |v(t)-v(s)|_{\FF L^{\alpha,p}}$ so that by dominated convergence  the function $u$ is only continuous $\FF L^{\alpha,p}$ without any further regularity. It is not difficult to prove that for $(\gamma,\alpha)$ in the interior of $\DD$ the solution $v$ is actually in $\CC^\gamma \FF L^{\alpha+\eps,p}$ for some small $\eps >0$. In this case it is clear that $u\in \CC^{\eps/3} \FF L^{\alpha,p}$ (cfr. the discussion of convergence of Galerkin approximation below and eq.~\eqref{eq:cancellation-breve-X2}).

Fix $p,\alpha$ and $\gamma$ such that $\gamma > 1/3$ and $(\gamma,\alpha) \in \DD'$. Then
an interesting property of the space $\QQ_\gamma$ is that for any continuous function $z$  in $\FF L^{\alpha,p}$ such that $y(t) = U(-t)z(t)$ is in 
$\QQ_\gamma$  the distribution $\mathcal{N}(P_N z)$ converge to a limit $\mathcal{N}(z)$. This follows from the proof of the previous corollary. Indeed  for general elements $y\in \QQ_\gamma$ the analog of eq.~\eqref{eq:weak-conv-decomp} reads
$$
\der V_{ts}=
\int_s^t U(-r)\mathcal{N}(P_N z(r)) dr
=Y_{ts}(y_s,y_s)+Y^2_{ts}(y_s,y'_s,y'_s) + V^\flat_{ts}
$$
and by the convergence of the couple $Y,Y^2$ to $X,X^2$ we have
$
\der V \to \der V^\infty
$ in $\CC^\gamma_1 \FF L^{\alpha,p}$ where $\der V^\infty$ is given by the $\Lambda$-equation $\der V^\infty = (1-\Lambda\der)[X(y,y)+X^2(y,y',y')]$. We can then define the distribution $\mathcal{N}(z)$ by letting $  \int_0^t U(-r)\mathcal{N}(z(r)) dr = V^\infty_{t}$ and have $\mathcal{N}(P_N z) \to \mathcal{N}(z)$ in weak sense.

\subsection{$L^2$ conservation law and global solutions}
\label{sec:global}
It is well known that the KdV equation formally conserves the $L^2(\TT)$ norm of the solution. This conservation law can be used to show
existence of global solution when the initial condition is in $L^2(\TT)$.  
Denote with $\langle \cdot,\cdot \rangle$ the $L^2$ scalar product.

\begin{lemma} 
\label{lemma:more-symmetry}
For $\varphi \in L^2$ we have $\langle \varphi, X_{ts}(\varphi,\varphi)\rangle=0$ and
$
2 \langle \varphi, X^2_{ts}(\varphi,\varphi,\varphi)\rangle
+   \langle X_{ts}(\varphi,\varphi), X_{ts}(\varphi,\varphi)\rangle = 0.
$  
\end{lemma}
\begin{proof} For smooth test functions 
we have
$\langle \varphi_1 , \dot X_s(\varphi_2,\varphi_2) \rangle  =
\int_{\TT} (U(s)\varphi_1)(\xi) (U(s)\varphi_2)(\xi) \partial_{\xi}(U(s)\varphi_2)(\xi) d\xi
$
and an integration by parts gives
$$
\langle \varphi_1 , \dot X_s(\varphi_2,\varphi_2) \rangle  =
- \int_{\TT} \partial_{\xi}[(U(s)\varphi_1)(\xi) (U(s)\varphi_2)(\xi)] (U(s)\varphi_2)(\xi) d\xi =
- 2 \langle \varphi_2 , \dot X_s(\varphi_1,\varphi_2) \rangle
$$
this gives directly that $\langle \varphi, X_{ts}(\varphi,\varphi)\rangle=0$. Moreover
  \begin{equation}
    \label{eq:cons-law-3}
    \begin{split}
  & \langle \varphi, X^2_{ts}(\varphi,\varphi,\varphi)\rangle = 
2 \int_s^t d\sigma \int_s^\sigma  d\sigma_1  
  \langle \varphi, \dot X_\sigma(\varphi,\dot
  X_{\sigma_1}(\varphi,\varphi))\rangle 
\\ & \qquad 
=- 
\int_s^t d\sigma \int_s^\sigma  d\sigma_1  
  \langle \dot
  X_{\sigma_1}(\varphi,\varphi), \dot X_\sigma(\varphi,\varphi)\rangle        
\\ & \qquad 
=- \frac12
\int_s^t d\sigma \int_s^t  d\sigma_1  
  \langle \dot
  X_{\sigma_1}(\varphi,\varphi), \dot X_\sigma(\varphi,\varphi)\rangle        
=- \frac12
  \langle 
  X_{ts}(\varphi,\varphi),  X_{ts}(\varphi,\varphi)\rangle        
    \end{split}
  \end{equation}
  and conclude by density.
\end{proof}

\begin{theorem}
If $v$ is a solution of eq.~(\ref{eq:abs-step-2}) in $[0,T_*]$ with initial condition $v_0 \in L^2$ then $|v_t|^2_{L^2(\TT)} = |v_0|^2_{L^2(\TT)}$ for any $t \in [0,T_*]$.
\end{theorem}
\begin{proof}
We will prove that $\der \langle v,v \rangle  = 0$. Let us compute explicitly this finite increment:
\begin{equation*}
  \begin{split}
[\der \langle v,v \rangle]_{ts}  & = \langle v_t,v_t \rangle - \langle v_s,v_s \rangle
=  2 \langle \der v_{ts},v_s \rangle + \langle \der v_{ts},\der v_{ts} \rangle    
  \end{split}
\end{equation*}
Substituting in this expression the $\Lambda$-equation~(\ref{eq:abs-step-2}) we get
\begin{equation*}
  \begin{split}
[\der \langle v,v \rangle]_{ts}
& = 2  \langle X_{ts}(v_s,v_s)+X^2_{ts}(v_s,v_s,v_s) + v^\flat_{ts},v_s \rangle 
\\ & \qquad + \langle X_{ts}(v_s,v_s),X_{ts}(v_s,v_s) \rangle 
+ 2 \langle X_{ts}(v_s,v_s),v^{\sharp}_{ts} \rangle    
+  \langle v^{\sharp}_{ts},v^{\sharp}_{ts} \rangle    
  \end{split}
\end{equation*}
where we set $v^\sharp = X^2(v,v,v) +
v^\flat$. Lemma~\ref{lemma:more-symmetry} implies that $\langle v_s,
X_{ts}(v_s,v_s) \rangle = 0$ and
allows to cancel the $X^2$ term with the quadratic $X$ term. After the
cancellations the increment of the $L^2$ norm squared is  
$[\der \langle v,v \rangle]_{ts}
 = 2  \langle  v^\flat_{ts},v_s \rangle 
+ 2 \langle X_{ts}(v_s,v_s),v^{\sharp}_{ts} \rangle    
+  \langle v^{\sharp}_{ts},v^{\sharp}_{ts} \rangle    
$. Each term  on the r.h.s. of this expression belongs at least to $\CC_2^{3\eta} \RR$ and
since $3\eta > 1$ this implies that the function $t \mapsto |v_t|^2_{L^2(\TT)}$
is an H\"older function of index greater than $1$ hence it
must be constant.
\end{proof}

\begin{corollary}
If $v_0 \in L^2$ there exist a unique global solutions to the $\Lambda$-equation~(\ref{eq:abs-step-2}).  
\end{corollary}
\begin{proof}
By Thm.~\ref{thm:main} there exists a unique local solution up to a time
$T_*$ which depends only on $|v_0|_{L^2}$. Since $|v_{T_*}|_{L^2} = |v_{0}|_{L^2}$
we can start from $T_*$ and extend uniquely this solution to the interval
$[0,2T_*]$ and then on any interval.
\end{proof}

It would be interesting to try to adapt the I-method
 of Colliander--Steel--Staffilani--Takaoka--Tao~\cite{MR1969209}  to extend the global well-posedness at least in the case $p=2$ for any $\alpha >  \alpha_*(2)= -1/2$. The handling of correction terms to the conservation law seems however to require some efforts and we prefer to leave this study to a further publication.

\subsection{Galerkin approximations}
\label{sec:galerkin}
Recall that $P_N$  is the projection on the Fourier modes $|k| \le N$. In~\cite{MR2233689} it is proven that the solutions of the approximate KdV equation
\begin{equation}
  \label{eq:kdv-real-N}
\partial_t u^{(N)} + \partial^3_\xi u^{(N)} +
\frac12 P_N \partial_\xi (u^{(N)})^2 = 0,
\qquad u^{(N)}(0) = P_N u_0
\end{equation}
do not converge even weakly to the flow of the full KdV equation. In the same paper the authors propose a modified finite dimensional scheme and prove its convergence in $H^{\alpha}(\TT)$ for any $\alpha \ge -1/2$. Here we would like to propose a different scheme inspired by the rough path analysis.
By partial series expansion for the twisted variable $v^{(N)}(t) = U(-t) u^{(N)}(t)$ is not difficult to show that the unique solution of equation~\eqref{eq:kdv-real-N} satisfy the $\Lambda$-equation
\begin{equation}
  \label{eq:abs-step-2-compact-N}
\der v^{(N)} = (\id - \Lambda \der)[X^{(N)}(v^{(N)})+X^{(N),2}(v^{(N)})]
\end{equation}
where 
$X^{(N)} = P_N X ( P_N \times P_N )$ and where the trilinear operator $X^{(N),2}$ is defined as
\begin{equation}
  \label{eq:X2-galerkin}
X^{(N),2}_{ts}(\varphi_1,\varphi_2,\varphi_3) = 2 \int_s^t d\sigma \int_s^\sigma d\sigma_1 P_N \dot X_\sigma(P_N \varphi_1,P_N \dot X_{\sigma_1}(P_N \varphi_2,P_N \varphi_3))  
\end{equation}
so that 
$
\der X^{(N),2}(\varphi_1,\varphi_2,\varphi_3) = 2 X^{(N)}(\varphi_1,X^{(N)}(\varphi_2,\varphi_3))
$.
These are just multi-linear operators in
a finite-dimensional space and to have convergence of the Galerkin
approximation it would be enough that both converge in norm to their
infinite-dimensional analogs $X,X^2$. 
A decomposition for   $X^{(N),2}$ analogous to that of $X^2$ described
in Lemma~\ref{lemma:X2} holds
$
X^{(N),2}_{ts} = \hat X^{(N),2}_{ts} + \breve X^{(N),2}_{ts}
$
and we will prove in the appendix that 
\begin{lemma}
\label{lemma:gal}
For any pair $(\gamma,\alpha)$ in the interior of $\DD$ we have that as  $N \to \infty$, $X^{(N)} \to X$ in $\CC_2^{\gamma} \LL^2 \FF L^{\alpha,p}$ and $\hat X^{(N),2} \to X^2$ in $\CC_2^{2\gamma} \LL^3 \FF L^{\alpha,p}$. 
\end{lemma}
Unfortunately it is not difficult to see that $\breve X^{(N),2}$ cannot converge in norm, indeed we have
$$
\FF
\breve X^{(N),2}_{ts}(\varphi_1, \varphi_{2}, \varphi_{3})(k) = (t-s) \sum_{k_1} \frac{ I_{0<|k|,|k_1|,|k_2| \le N}}{3 i k_1} \wh \varphi_1 (k_1)[\wh \varphi_{2} (-k_1) \wh \varphi_{3} (k)+\wh \varphi_{3} (-k_1) \wh \varphi_{2} (k)]
$$
and there is no way to make this converge in norm to $\breve X_{ts}^2=\breve X_{ts}^{(\infty),2}$ due to the explicit dependence of the cutoff on $|k|,|k_2|$ which cannot be compensated by the regularity of the test functions or by the $1/k_1$ factor.
A way to remove this difficulty is to modify the finite dimensional ODE in order to remove this operator in the $\Lambda$-equation. This is possible since $\breve X^{(N),2}_{ts} $ is proportional to $t-s$ so that it admits an obvious differential counterpart. Let use define the trilinear operator $\Gamma^{(N)}$ as
$$
\FF \Gamma^{(N)} (\varphi_1, \varphi_{2}, \varphi_3)(k) =  \sum_{k_1} \frac{(I_{0<|k|,|k_1|,|k_2|}- I_{0<|k|,|k_1|,|k_2| \le N})}{3 i k_1} \wh \varphi_1 (k_1)[\wh \varphi_{2} (-k_1) \wh \varphi_{3} (k)+\wh \varphi_{3} (-k_1) \wh \varphi_{2} (k)]
$$
and note that 
$$
\breve X^{2}_{ts}(\varphi_1, \varphi_{2}, \varphi_{3})-\breve X^{(N),2}_{ts}(\varphi_1, \varphi_{2}, \varphi_{3})  = \int_s^t U(-r)\Gamma(U(r)\varphi_1, U(r)\varphi_{2}, U(r)\varphi_{3}) dr .
$$
Then the modified Galerkin scheme 
\begin{equation}
  \label{eq:kdv-real-N-mod}
\partial_t u^{(N)} + \partial^3_\xi u^{(N)} +
\frac12 P_N \partial_\xi (u^{(N)})^2 - \Gamma^{(N)}(u^{(N)})= 0,
\qquad u^{(N)}(0) = P_N u_0
\end{equation}
is still finite dimensional since $\Gamma^{(N)}P_N^{\times 3} = P_N \Gamma^{(N)}P_N^{\times 3}$ and
is equivalent to the $\Lambda$-equation
\begin{equation}
\begin{split}
  \label{eq:abs-step-2-compact-N-bis}
\der v^{(N)} & = (\id - \Lambda \der)[X^{(N)}(v^{(N)})+X^{(N),2}(v^{(N)})
 -\breve X^{(N),2}(v^{(N)})+\breve X^{2}(v^{(N)})]
\\ &= (\id - \Lambda \der)[X^{(N)}(v^{(N)})+\tilde X^{(N),2}(v^{(N)})]
\end{split}
\end{equation}
where $\tilde X^{(N),2} = \hat X^{(N),2} + \breve X^{2}$ now do converge in norm to $X^2$ and satisfy the correct algebraic relations.
As a consequence of the Lipschitz continuity of the solution of the
$\Lambda$-equation~(\ref{eq:abs-step-2-compact-N-bis}) w.r.t $X, X^2$ and $v_0$ implies
the following convergence result.

\begin{corollary}
\label{cor:galerkin} Let $1 \le p \le +\infty$ and $\alpha > \alpha_*(p)$, then for any $u_0 \in \FF L^{\alpha,p}$
as $N \to \infty$ the Galerkin approximations $v^{(N)}_t = U(-t) u^{(N)}_t$ obtained by the ODE~\eqref{eq:kdv-real-N-mod}  converges in $\CC_1^\gamma  \FF L^{\alpha,p}$ to
the solution $v$ of the $\Lambda$-equation~(\ref{eq:abs-step-2}) up to
a strictly positive time $T^*$ which depends only on the norm of
$v_0,X,X^2$.  
\end{corollary}
It is whortwhile to note that this result imply that
$
\sup_{t\in [0,T^*]} |u^{(N)}(t) - u(t)|_{\FF L^{\alpha,p}}  \to 0
$
as $N \to \infty$ while for the finite dimensional scheme devised in~\cite{MR2233689} the convergence holds only in the sense that
$
\sup_{t\in [0,T^*]} |P_{\sqrt N }(u^{(N)}(t) - u(t))|_{H^{\alpha}}  \to 0
$
i.e. only for a very low frequency part of the solution.
It is interesting to remark that the modified ODE~\eqref{eq:kdv-real-N-mod} remains an Hamiltonian flow on $P_N (H^{-1/2}(\TT)\backslash \RR)$ endowed with the
symplectic structure given by
$
\Omega(u,v) = \sum_{0 < |k| \le N} u(-k) v(k)/(ik)
$. Its Hamiltonian is given by
$$
H(u) = \frac 1 2 \sum_{0< |k| \le N}k^2 u(-k) u(k) + \frac 1 6\sum_{\substack{k_1+k_2+k_3=0\\ 0 < |k_i| \le N}} u(k_1) u(k_2) u(k_3)
$$
$$\qquad - \frac 1 {12} \sum_{0< |k|,|k_1| \le N}  \frac{u(-k) u(k)}{ik} \frac{u(-k_1) u(k_1)}{i k_1}  I_{|k-k_1|>N} .
$$

\subsection{A discrete time scheme}
\label{sec:euler}
The solution described by the
$\Lambda$-equation~(\ref{eq:abs-step-2}) can be approximated
by a discrete Euler-like scheme defined as follows. For any $n >0$ let $y^n_0=v_0$
and 
$$
y^n_i = X_{i/n, j/n}(y^n_{i-1}) +  X^2_{i/n, j/n}(y^n_{i-1})
$$
for $i\ge 1$. 
The combination of  this scheme with the Galerkin approximation
discussed before provide an implementable numerical  approximation scheme for the
solutions of KdV with low regularity initial conditions. Indeed the next theorem  can  be combined with
Corollary~\ref{cor:galerkin} to obtain effective rates of convergence.

\begin{theorem}
Let $\Delta^n_i = y^n_i-y_{i/n}$ and let $T>0$ be the existence time of the
solution described in Thm.~\ref{thm:main}, then
\begin{equation}
  \label{eq:euler-bound}
\sup_{0\le i< j\le nT} \frac{|\Delta^n_i-\Delta^n_j|_{\FF L^{\alpha,p}}}{|i-j|^\gamma} = O(n^{1-3\gamma}) . 
\end{equation}
\end{theorem}

\begin{proof}
Let $T$ be the existence time for the solution $v$ described in theorem~\ref{thm:main} and let $N = \lfloor nT \rfloor$.
We begin by proving some uniform bounds on the sequence
$\{y^n_i,i=0,\dots,N\}$.   Let
$$
q^n_{ij} = X_{i/n, j/n}(y^n_{i}) +  X^2_{i/n, j/n}(y^n_{i})
$$
so that $y^n_k =\sum_{i=0}^{k-1} q_{i,i+1}$. Given $0\le i<j\le N$ let
$\tau^1_0=i, \tau^1_1=j$ and define recursively  $\tau^k_l$, $k >
0$, $l =0,\dots,2^k$ such that $\tau^{k+1}_{2l} = \tau^k_l$  for
$l=0,\dots,2^k$ and 
$\tau^{k+1}_{2l+1}<\lfloor (\tau^{k}_{l}+\tau^{k}_{l+1})/2\rfloor$ for
$l=0,\dots,2^k-1$. 
Then we have $|\tau^k_{l+1}-\tau^k_l|\le 1 \vee (j-i)/2^k$
and
$
y^n_j-y^n_i =\sum_{l=0}^{2^K} q^n_{\tau^k_l,\tau^k_{l+1}} 
$
 where $K$ is such that
$(j-i)/2^K \le 1$. Using the triangular array $\tau^k_l$ we 
rewrite the above expression as a telescopic sum:
\begin{equation*}
  \begin{split}
y^n_j-y^n_i & = q^n_{ij} + \sum_{k=1}^K \sum_{l=0}^{2^{K-1}} (q^n_{\tau^k_{2l},\tau^k_{2l+2}}-q^n_{\tau^k_{2l},\tau^k_{2l+1}}-q^n_{\tau^k_{2l+1},\tau^k_{2l+2}}) 
  \\ & = q^n_{ij} + \sum_{k=1}^K \sum_{l=0}^{2^{k-1}} (\der q^n)_{\tau^k_{2l},\tau^k_{2l+1},\tau^k_{2l+2}}   
  \end{split}
\end{equation*}
Up to $T$, the solution $v$ satisfy the equation
$
\der v = X(v) + X^2(v) + v^\flat
$
where $\sup_{t\in[0,T]}|v_t|_\alpha+\|v^\flat\|_{3\gamma}\le C$ so
$
v_{j/n}-v_{i/n} = \sum_{l=i}^{j-1} p_{l,l+1} + \sum_{l=i}^{j-1} v^\flat_{l/n,(l+1)/n}
$
where $p_{ij}=X_{i/n,j/n}(v_{i/n})+ X^2_{i/n,j/n}(v_{i/n})$.
Then
\begin{equation*}
\Delta^n_{j}-\Delta^n_{i} =  q^n_{ij}-p_{ij}  + \sum_{k=1}^K \sum_{l=0}^{2^{k-1}} (\der
q^n- \der p)_{\tau^k_{2l},\tau^k_{2l+1},\tau^k_{2l+2}}  - r_{ij}
\end{equation*}
where
$
r_{ij}=  \sum_{l=i}^{j-1} v^\flat_{l/n,(l+1)/n}
$.
This last term is readily estimated by
$$
|r_{ij}| \le \sum_{l=i}^{j-1} \|v^\flat\|_{3\gamma} n^{-3\gamma}
 \le C \left(\frac{j-i}{ n}\right) n^{1-3\gamma}
$$
uniformly in $n$.
Let 
$$
M^n_\ell = \sup_{0\le i < j \le \ell} \left(\frac{j-i}{ n}\right)^{-1} |\Delta^n_{j}-\Delta^n_{i} -  q^n_{ij}+p_{ij}+r_{ij}|
$$
we want to show that $M^n_N \le A n^{1-3\gamma}$ uniformly in $n$ for some constant
$A$ depending only on the data of the problem: this will imply the
statement of the theorem since then
$$
|\Delta^n_{i}| \le |q^n_{0i}-p_{0i}| + |r_{0i}| + A n^{1-3\gamma}
\left(\frac{i}{ n}\right)  \le |r_{0i}| + A n^{1-3\gamma}
\left(\frac{i}{ n}\right) \le Cn^{1-3\gamma}
$$
for any $i \le N$ and
$$
|\Delta^n_{j}-\Delta^n_{i}| \le |q^n_{ij}-p^{ij}| + |r_{ij}| + A n^{1-3\gamma}
\left(\frac{j-i}{ n}\right)  \le C n^{1-3\gamma} \left(\frac{j-i}{ n}\right)^\gamma
$$
for any $0\le i<j\le N$.

 We proceed by
induction on $\ell$. For $\ell=1$ the statement is clearly true since
$\Delta^n_{1}-\Delta^n_{0} -  q^n_{0,1}+p_{0,1}+r_{0,1}=0$, moreover
for the same reason we have, for all $l$ that
$\Delta^n_{l+1}-\Delta^n_{l} -  q^n_{l,l+1}+p_{l,l+1}+r_{l,l+1}=0$.
 Assume then
that $M^n_{\ell-1} \le A$ for some $\ell> 0$. The basic observation is that when $|i-j|\le \ell$
the sums
$$
\Delta^n_{j}-\Delta^n_{i} -  q^n_{ij}+p_{ij}+r_{ij} = \sum_{k=1}^K \sum_{l=0}^{2^{k-1}} (\der
q^n- \der p)_{\tau^k_{2l},\tau^k_{2l+1},\tau^k_{2l+2}}
$$
can be estimated in terms of $M^n_{\ell-1}$ and various norms of
$v$,$X$,$X^2$ much like in the proof of Thm.~\ref{thm:main}. The bound
has the form
$$
|(\der
q^n- \der p)_{\tau^k_{2l},\tau^k_{2l+1},\tau^k_{2l+2}}| \le C
\left(\frac{j-i}{2^k n}\right)^{3\gamma} (1+M^n_{\ell-1})^3 (M^n_{\ell-1}+n^{1-3\gamma})
$$
where $C=C(v_0,X,X^2)$ and where the factor $n^{1-3\gamma}$ is due to
the previous estimate on $r_{ij}$. Then since $3\gamma > 1$ and
$$
\sum_{k=1}^K \sum_{l=0}^{2^{k-1}}
\left(\frac{j-i}{2^k n}\right)^{3\gamma}  \le \left(\frac{j-i}{
    n}\right)^{3\gamma} \sum_{k=1}^\infty 2^{k(1-3\gamma)} \le C  \left(\frac{j-i}{
    n}\right)^{3\gamma} \le C \left(\frac{\ell}{
    n}\right)^{3\gamma-1}  \left(\frac{j-i}{
    n}\right) 
$$
we get
$$
M^n_\ell \le C
(1+M^n_{\ell-1})^4 \sum_{k=1}^K \sum_{l=0}^{2^{k-1}}
\left(\frac{j-i}{2^k n}\right)^{3\gamma} \le C (1+M^n_{\ell-1})^3 (M^n_{\ell-1}+n^{1-3\gamma}) \left(\frac{\ell}{
    n}\right)^{3\gamma-1} 
$$
Let $m^n_\ell = n^{3\gamma-1} M^{n}_\ell$, then, for $n$ large enough
$m^n_\ell \le F(m^n_{\ell-1})$ where
$F$ is the increasing map
$$
\RR_+ \ni m \mapsto F(m)= C (1+m)^4 \left(\frac{\ell}{
    n}\right)^{3\gamma-1} \in \RR_+
$$
which, for $\ell/n$ small enough
has a unique attracting fix-point under iteration starting from $0$. 
 In particular the
iterations stay bounded and if we  set $x_{i+1}=F(x_i)$, $x_0=0$ we
have $A = \sup_i F(x_i)  <
\infty$ and $m^n_\ell \le x_\ell \le A$. By repeating this argument it is easy to prove that the
bound holds for all $\ell \le nT$, i.e. in the whole existence
interval found in Thm.~\ref{thm:main}.
\end{proof}

\begin{remark}
With a bit more of work it is possible to prove the existence of
the solution stated in Thm.~\ref{thm:main} using directly the discrete
approximation as done by Davie~\cite{davie} for rough
differential equations.
\end{remark}

\subsection{Higher order $\Lambda$ equations}
\label{sec:more}
Further expansion of eq.~(\ref{eq:exp-step-2}) generate a hierarchy of
$\Lambda$-equations for KdV. The next one is given by
\begin{equation}
  \label{eq:hier-next}
  \der v = (\id-\Lambda\der) [X(v) + X^2 (v) + X^{3a}(v) + X^{3b}(v)]
\end{equation}
 where $X^{3a},X^{3b}$ are operators increments respectively defined as
\begin{equation}
  \label{eq:X3a}
  X^{3a}_{ts}(\varphi_1,\varphi_2,\varphi_3,\varphi_4) = \int_s^t
  d\sigma_1 \dot
  X_{\sigma_1}(\varphi_1,\int_s^{\sigma_1}d\sigma_2 \dot X_{\sigma_2}(\varphi_2,
\int_s^{\sigma_2}d\sigma_3 \dot X_{\sigma_3}(\varphi_3,\varphi_4)))
\end{equation}
and
\begin{equation}
  \label{eq:X3b}
  X^{3b}_{ts}(\varphi_1,\varphi_2,\varphi_3,\varphi_4) = \int_s^t
  d\sigma_1 \dot
  X_{\sigma_1}(\int_s^{\sigma_1}d\sigma_2 \dot
  X_{\sigma_2}(\varphi_1,\varphi_2),
\int_s^{\sigma_1}d\sigma_3 \dot
  X_{\sigma_3}(\varphi_3,\varphi_4))
\end{equation}
which satisfy the following relations with $X$ and $X^2$:
\begin{equation}
  \label{eq:deltaX3a}
  \der X^{3a}(\varphi_1,\varphi_2,\varphi_3,\varphi_4) =
  X(\varphi_1,X^2(\varphi_2,\varphi_3,\varphi_4)) + X^2(\varphi_1,\varphi_2,X(\varphi_3,\varphi_4))
\end{equation}
and
\begin{equation}
  \label{eq:deltaX3b}
  \begin{split}
  \der X^{3b}(\varphi_1,\varphi_2,\varphi_3,\varphi_4) & =
  X(X(\varphi_1,\varphi_2),X(\varphi_3,\varphi_4)) +
  X^2(\varphi_3,\varphi_4,X(\varphi_1,\varphi_2))
\\ & \qquad  + X^2(\varphi_1,\varphi_2,X(\varphi_3,\varphi_4))    
  \end{split}
\end{equation}
As we report elsewhere~\cite{trees}  the Hopf algebra of rooted trees is
the natural language to describe this hierarchy of equations and the
algebraic relations between the various operators. In this special
case however these relations can be easily checked by direct 
 computations. 
Using Lemmas~\ref{lemma:X} and~\ref{lemma:X2} we can show that the r.h.s. of
the eqns.~(\ref{eq:deltaX3a}) and~(\ref{eq:deltaX3b}) belongs to the
domain of $\Lambda$ and so the equations can be used to express
$X^{3a},X^{3b}$ in function of $X,X^2$ and prove that
\begin{corollary}
For $(\gamma,\alpha) \in\DD \cap \{\alpha \ge -1/p\}$ with $\gamma>1/3$ we have
$X^{3a},X^{3b} \in \CC_2^\gamma \LL^4 \FF L^{\alpha,p}.$ 
\end{corollary}

\section{Additive stochastic forcing}
\label{sec:stoch}
As another application of this approach we would like to discuss the presence of an additive random force in the KdV eq.~(\ref{eq:kdv-real}):
\begin{equation}
  \label{eq:kdv-real-stoch}
\partial_t u + \partial^3_\xi u +
\frac12 \partial_\xi u^2 = \Phi \partial_t  \partial_\xi B
\end{equation}
where 
$\partial_t \partial_\xi B$ a white noise on $\RR\times\TT$ and where
$\Phi$ is a linear operator acting on the space variable which is diagonal in Fourier space: $\Phi e_k =
\lambda_k e_k$ where $\{e_k \}_{k\in\ZZ}$ is the orthonormal basis
$e_k(\xi) = e^{ik\xi}/\sqrt{2\pi}$ and such that $\lambda_0 = 0$.  In this
way the noise does not affect the zero mode.
In the rest of this section  fix $1\le p \le +\infty$, $\alpha > \alpha_*(p)$ and $\gamma\in (1/3,1/2)$ such that $(\gamma,\alpha)
\in \DD'$ and assume that
\begin{equation}
  \label{eq:lambda-bound}
|\lambda|_{\ell^{\alpha,p}} = \sum_{k \in \ZZ}  |k|^{\alpha p} |\lambda_k|^p < \infty .
\end{equation}
The transformed integral equation (analogous to
eq.~(\ref{eq:kdv-base})) associated to~(\ref{eq:kdv-real-stoch}) is
\begin{equation}
  \label{eq:kdv-base-stoch}
 v_t =   v_0 + w_t +  \int_0^t \dot X_s(
  v_s, v_s) \, ds
\end{equation}
where $\wh w_t(k) = \lambda_k \beta^k_t$ and $\{\beta^k_\cdot\}_{k \in \ZZ_*}$ is a family of 
complex-valued centered Brownian motions such that $\beta^{-k} =
\overline{\beta^k}$ and with covariance
$\expect[\overline \beta^k_t \beta^q_s] = \delta_{k,q} (t\wedge s)$. 
The relation between the initial noise $B$ and the family $\{\beta^k_\cdot\}_{k \in \ZZ_*}$ is given by 
$
\beta^k_t = \langle e_k, \int_0^t U(-s) \partial_\xi \partial_s B(s,\cdot) ds \rangle_{H^0} 
$.
Eq.~(\ref{eq:kdv-base-stoch}) can be expanded in the same
way as we have done before and the first
interesting $\Lambda$-equation which appears is the following:
\begin{equation}
  \label{eq:lambda-eq-stoch}
  \der v = (\id-\Lambda\der)[X(v) + \der w + X^{2}(v) + X^{w}(v)].
\end{equation}
Here   the
random operator $X^w_{ts}: \FF L^{\alpha,p}\to \FF L^{\alpha,p}$ is given by
\begin{equation}
  \label{eq:Xw}
  X^w_{ts}(\varphi) = \int_s^t d\sigma \dot X_\sigma(\varphi, \der
  w_{\sigma s})
\end{equation}
and satisfy the equation
$  \der X^w_{tus} = X_{tu}(\varphi,\der w_{us})
$.
For any couple of integers $n,m$ we have
\begin{equation}
\label{eq:w-norm-est}
\begin{split}
\expect & |\der w_{ts}|_{\FF L^{\alpha,2n}}^{2m n}  = \expect \left[ \sum_k |k|^{\alpha 2n} |\der \wh w_{ts}(k)|^{2n} \right]^{m} 
= \sum_{k_1,..,k_m}  \expect \left[\prod_{i=1}^m |k_i|^{\alpha 2n} |\der \wh w_{ts}(k_i)|^{2n} \right]
\\
& \le \sum_{k_1,..,k_m}  \prod_{i=1}^m   \left[|k_i|^{\alpha 2nm}\expect |\der \wh w_{ts}(k_i)|^{2nm} \right]^{1/m}
 = \left\{\sum_{k}   \left[|k|^{\alpha 2nm}\expect |\der \wh w_{ts}(k)|^{2nm} \right]^{1/m}\right\}^m
\\
& \lesssim_{nm} \left\{\sum_{k}  |k|^{\alpha 2n} |\lambda_k|^{2n}\right\}^m |t-s|^{nm} = \|\lambda\|_{\ell^{\alpha,2n}}^{2nm} |t-s|^{nm}
\end{split}
\end{equation}
where we used the Gaussian bound
$$
\expect |\der \wh w_{ts}(k)|^{2nm} \le C_{nm} \left(\expect |\der \wh w_{ts}(k)|^{2}\right)^{nm} \le C_{nm} |\lambda_k|^{2nm} .
$$
By interpolation this gives $\expect  |\der w_{ts}|_{\FF L^{\alpha,p}}^{r} \lesssim  |\lambda|_{\ell^{\alpha,p}}^{r} |t-s|^{r/2}$ for all $r\ge p \ge 2$ which is finite by assumption~(\ref{eq:lambda-bound}).
By the standard Kolmogorov criterion this implies that a.s. $\der w \in
\CC^{\rho}_2 \FF L^{\alpha,p}$ for any $\rho < 1/2$ and a-fortiori $\der w \in
\CC^{\gamma}_2 \FF L^{\alpha,p}$ by choosing $\rho \in [\gamma,1/2)$. 
To prove that a sufficiently regular version of the Gaussian stochastic process $X^w$ exists
we will use a generalization of the classic Garsia-Rodemich-Rumsey lemma
which has been proved in~\cite{MR2091358}.
\begin{lemma}
\label{lemma:besov}
For any $\theta > 0$ and $p \ge 1$,  there exists a constant
$C$ such that for any $R \in \CC_{2} V$ ($(V,|\cdot|)$ some Banach space), we have 
\begin{equation}
\label{eq:generalboundxx}
\|R\|_{\theta} \le C (U_{\theta+2/p,p}(R)+\|\der R\|_{\theta}),
\end{equation}
where
$
 U_{\theta,p}(R) = \left[ \iint_{[0,T]^2}
 \left(\frac{|R_{t s}|}{|t-s|^\theta}\right)^p dt ds \right]^{1/p}.
$
\end{lemma}

The operator $X^w$ behaves not worse than $X^2$: 

\begin{lemma}
\label{lemma:Xw}
Under condition~(\ref{eq:lambda-bound})
we have $X^w \in \CC^{2\gamma}_2 \LL \FF L^{\alpha,p}$ a.s..
\end{lemma}
\begin{proof}
After an integration by parts, $X^w$ can be rewritten as
\begin{equation}
  \label{eq:Xw-2}
  X^w_{ts}(\varphi) = X_{ts}(\varphi,\der w_{ts}) - \int_s^t X_{\sigma
  s}(\varphi,dw_\sigma) = X_{ts}(\varphi,\der w_{ts}) + I_{ts}(\varphi) 
\end{equation}
The first term in the r.h.s. belongs to $\CC_2^{2\gamma}
\LL \FF L^{\alpha,p}$ path-wise:
$$
|X_{ts}(\varphi,\der w_{ts})|_{\FF L^{\alpha,p}} \le
|X_{ts}|_{\LL \FF L^{\alpha,p}} |\der w_{ts}|_{\alpha,p}
|\varphi|_{\alpha,p} \le 
\|X\|_{\CC^\gamma_2 \LL \FF L^{\alpha,p}} \|\der
w\|_{\CC_1^\gamma \FF L^{\alpha,p}} |\varphi|_{\FF L^{\alpha,p}} |t-s|^{2\gamma}.
$$
Let us estimate the random operator $I_{ts} : \varphi \mapsto I_{ts}(\varphi)$. Its Fourier kernel is
$$
\FF I_{ts}(\varphi)(k) = \int_{s}^t \sum_{k_1} \frac{e^{-i 3kk_1 k_2 \sigma}-e^{-i 3k k_1 k_2 s}}{6 k_1 k_2}  \lambda_{k_2} \wh\varphi(k_1) d\beta^{k_2}_\sigma
=  \sum_{k_1} \frac{|k k_1 k_2|^\gamma}{6 k_1 k_2} \lambda_{k_2} \wh\varphi(k_1) J(k,k_1,k_2)
$$
where 
$$
J_{ts}(k,k_1,k_2) =  \int_{s}^t \frac{e^{-i 3kk_1 k_2 \sigma}-e^{-i 3k k_1 k_2 s}}{ |k k_1 k_2|^{\gamma}}   d\beta^{k_2}_\sigma
$$
so
$$
|I_{ts}|_{\LL \FF L^{\alpha,p}} \le |\YY_2(Q)|_{\FF L^p} |\lambda|_{\ell^{\alpha,p}} \sup_{k,k_1} |J_{ts}(k,k_1,k_2)|
$$
with 
$$
Q(k,k_1,k_2) = \frac{|k|^{\alpha+\gamma}}{|k_1 k_2|^{1+\alpha-\gamma}} I_{k_1 k_2 \neq 0} .
$$
The majorizing kernel $Q$ is the same appearing in the estimates for $X$ so that we already know that $\|\YY_2(Q)\|_{\FF L^p} < \infty$ for all allowed pairs $(\gamma,\alpha)$. It remains to show that
$
\expect \sup_{k,k_1} |J_{ts}(k,k_1,k_2)|^n \lesssim |t-s|^{n(1+\gamma)} 
$
for arbitrarily large $n$. It is then enough to bound
\begin{equation*}
\begin{split}
\expect \sup_{k,k_1} |J_{ts}(k,k_1,k_2)|^{2n} & \le \sum_{k,k_1} \expect  |J_{ts}(k,k_1,k_2)|^{2n}
\lesssim_n \sum_{k,k_1} \left[\int_s^t \frac{|e^{-i 3kk_1 k_2 \sigma}-e^{-i 3k k_1 k_2 s}|^2}{ |k k_1 k_2|^{2\gamma}} d\sigma  \right]^n
\\†& \lesssim_n |t-s|^{n+2 \gamma'} \sum_{k_1,k_2}  \frac{1}{ |k k_1 k_2|^{2n(\gamma-\gamma')}} 
\end{split}
\end{equation*}
where the sum is finite for $n$ large enough (depending on $\gamma-\gamma'$).
Then choosing $p$ sufficiently large, we
have $U_{2\gamma+2/p,p}(I) < \infty$ a.s.. Moreover $\der I_{tus}(\varphi)
= X_{tu}(\varphi,\der w_{us})$ so that $\der I$ can be bounded
path-wise in $\CC_3^{2\gamma} \LL \FF L^{\alpha,p}$.  
Then Lemma~\ref{lemma:besov} implies that 
$
\|I\|_{2\gamma} \le C(U_{2\gamma+2/p,p}(I)+\|\der I\|_{2\gamma}) < \infty
$ a.s.
ending the proof.
\end{proof}

Then, modifying a bit the proof of Thm.~\ref{thm:main}, is not
difficult to prove the following.

\begin{theorem}
\label{th:eq-stoch} For any $1 \le p \le +\infty$ and $\alpha > \alpha_*(p)$
eq.~(\ref{eq:lambda-eq-stoch}) has a unique local path-wise solution in $\CC^\gamma \FF L^{\alpha,p}$ for any initial
condition in $\FF L^{\alpha,p}$.  
\end{theorem}
When $p=2$  we obtain
solutions for noises with values in $H^\alpha(\TT)$ for any  $\alpha > -1/2$.
In this way we essentially cover and extend the results of
De~Bouard-Debussche-Tsutsumi~\cite{MR2111917}. Their approach consist
in modifying Bourgain's method to handle  Besov spaces in order to compensate
for the insufficient  Sobolev time regularity of Brownian motion.

\section*{Acknowledgment}
I would like to thank A. Debussche which  delivered a series of interesting
lectures on stochastic dispersive equations  during a 2006
semester on Stochastic Analysis at Centro de Giorgi, Pisa. They
constituted the motivation for the investigations reported in this
note. I'm also greatly indebted with J.~Colliander and with an anonymous referee for  some remarks which helped me to discover an error in an earlier version of the paper. 


\appendix
\section{Regularity of some operators}
\label{sec:reg}

Some elementary results needed in the proofs of this appendix are the
subject of the next few lemmas.
We have to deal with $n$-multilinear operators $\YY_n(m) : (\FF L^p)^{ n}\to \FF L^p$ associated to multipliers $m:\RR^n\to\mathbb{C}$ as
$$
\FF [\YY_n(m)(\psi_1,\dots,\psi_n)](k_0) = \sum_{k_0+k_1+\cdots+k_n=0} m(k_0,k_1,\dots,k_n) \hat \psi_1(k_1)\cdots \hat \psi_n(k_n) .
$$
Recall that we have the interpolation inequalities
$$
|\YY_n(m(t))|_{\LL^n \FF L^{\alpha,p(t)}} \le |\YY_n(m_1) |_{\LL^n \FF L^{\alpha,p_1}}^{t} |\YY_n(m_2)|_{\LL^n \FF L^{\alpha,p_2}}^{1-t}
$$
for any positive multiplier $m(t) = m_1^t m_2^{1-t}$ and $t\in[0,1]$ such that  $\YY_n(m_1) \in \LL^n \FF L^{\alpha,1}$ and $\YY_n(m_\infty) \in \LL^n\FF L^{\alpha,\infty}$ and where $1/p(t) = t/p_1+(1-t)/p_2$. 
The $\YY$ operators can be bounded   in $\LL^n \FF L^p$ in terms of the multipliers as stipulated  by the following
\begin{lemma}
\label{lemma:op-bound}
For $1/p+1/q=1$ we have
$$
|\YY_n(m)|_{\LL^n\FF L^p}\le \sup_{k_0} \big[\sum_{\substack{k_0+k_1+\cdots+k_n=0\\  k_0 \text{ fixed}}} |m(k_0,k_1,\dots,k_n)|^q \big]^{1/q}
$$
and, for any $p\ge n/(n-1)$,
$$
|\YY_n(m)|_{\LL^n\FF L^p}
\le \big 
\{\sum_{k_0}\big[\sum_{\substack{k_0+k_1+\cdots+k_n=0\\ k_0 \text{ fixed}}} |m|^{\hat q} \big]^{p/\hat q} \big\}^{1/p}
$$
where $\hat q = p/[p-n/(n-1)]$.

\end{lemma}
\begin{proof}
By duality it is enough to bound the linear form
$$
F={}_{\FF L^q}\langle \psi_0,\YY_n(m)(\psi_1,\dots,\psi_n)\rangle_{\FF L^p}= \sum_{k_0+k_1+\cdots+k_n=0} m(k_0,k_1,\dots,k_n) \hat \psi_0(k_0)\cdots \hat \psi_n(k_n)
$$
for any $\psi_0 \in \FF L^q$:
$$
|F|\le \sum_{k_0}| \hat \psi_0| \big[\sum_{\substack{k_0+k_1+\cdots+k_n=0\\ k_0 \text{ fixed}}} |m |^q \big]^{1/q}\big[\sum_{\substack{k_0+k_1+\cdots+k_n=0\\ k_0 \text{ fixed}}} |\hat \psi_1 \cdots \hat \psi_n|^p\big]^{1/p}
$$
$$
\le \sup_{k_0} \big[\sum_{\substack{k_0+k_1+\cdots+k_n=0\\ k_0 \text{ fixed}}} |m |^q \big]^{1/q} \sum_{k_0} |\hat \psi_0|
\big[\sum_{\substack{k_0+k_1+\cdots+k_n=0\\ k_0 \text{ fixed}}} | \hat \psi_1\cdots \hat \psi_n|^p\big]^{1/p}
$$
$$
\le \sup_{k_0} \big[\sum_{\substack{k_0+k_1+\cdots+k_n=0\\ k_0 \text{ fixed}}} |m |^q \big]^{1/q}  |\psi_0|_{\FF L^q}  |\psi_1|_{\FF L^p}\cdots  |\psi_n|_{\FF L^p}
$$
For the second inequality we have
$$
|F|\le\sum_{k_0}  |\hat \psi_0|  \sum_{\substack{k_0+k_1+\cdots+k_n=0\\ k_0 \text{ fixed}}} |m| \prod_{k=1}^n |\hat \psi_k|
$$
$$
\le
\sum_{k_0}  |\hat \psi_0| [\sum_{\substack{k_0+k_1+\cdots+k_n=0\\ k_0 \text{ fixed}}} |m|^{\hat q}]^{1/\hat q}[\sum_{\substack{k_0+k_1+\cdots+k_n=0\\ k_0 \text{ fixed}}} \prod_{a=1}^n |\hat \psi_a|^{\hat p}]^{1/\hat p}
$$
Now choosing $n \hat p = (n-1)p$ we have
$$
\sum_{\substack{k_0+k_1+\cdots+k_n=0\\ k_0 \text{ fixed}}} \prod_{k=1}^n |\hat \psi_k|^{\hat p}
=
\sum_{\substack{k_0+k_1+\cdots+k_n=0\\ k_0 \text{ fixed}}} \prod_{j=1}^n(\prod_{\substack{a=1,..,n\\a\neq j}} |\hat \psi_a|^{\hat p/(n-1)})
$$
$$
\le\prod_{j=1}^n \big[ \sum_{\substack{k_0+k_1+\cdots+k_n=0\\ k_0 \text{ fixed}}} \big(\prod_{\substack{k=1,..,n\\k\neq j}} |\hat \psi_k|^{n \hat p/(n-1)}\big)\big]^{1/n}
$$
$$
=\prod_{j=1}^n \big[ \sum_{\substack{k_1,\dots,\hat k_j,\dots,k_n}} \big(\prod_{\substack{a=1,..,n\\a\neq j}} |\hat \psi_a|^{n \hat p/(n-1)}\big)\big]^{1/n}=\prod_{j=1}^n |\psi_j |_{\FF L^p}^{\hat p}
$$
so that
$$
|F|\le\prod_{j=1}^n |\psi_j |_{\FF L^p}
\sum_{k_0}  |\hat \psi_0|  \big[\sum_{\substack{k_0+k_1+\cdots+k_n=0\\ k_0 \text{ fixed}}} |m|^{\hat q} \big]^{1/\hat q}
$$
$$ \le \prod_{j=1}^n |\psi_j |_{\FF L^p} |\psi_0|_{\FF L^q} \big 
\{\sum_{k_0}\big[\sum_{\substack{k_0+k_1+\cdots+k_n=0\\ k_0 \text{ fixed}}} |m|^{\hat q} \big]^{p/\hat q} \big\}^{1/p}
$$
with $\hat q = p/[p-n/(n-1)]$.
\end{proof}

\begin{lemma}
\label{lemma:M-bounds}
 Fix any $a,b \in \RR$ and let 
$
m = \int_s^t \int_s^\sigma e^{i a \sigma} e^{i b \sigma_1} d\sigma d\sigma_1 - I_{a+b=0} (t-s)/(ib)
$.
If $\gamma \in[0,1/2]$  we have
$$
|m| \lesssim\frac{|t-s|^{2\gamma}}{|b|^{\gamma}|a|^{1-\gamma}|a+b|^{1-2\gamma}}+\frac{|t-s|^{2\gamma}}{|a|^{1-\gamma}|b|^{1-\gamma}} .
$$  
\end{lemma}
\begin{proof}
An explicit integration gives
$$
m = \frac{e^{i(a+b)s}-e^{i(a+b)t}}{b (a+b)} 
+\frac{e^{i bs+iat}-e^{i(a+b)s}}{ab}  - I_{a+b=0} \frac{t-s}{ib}.
$$
Then if $a+b \neq 0$
\begin{equation}
  \label{eq:m-bound-1}
|m| \le \frac{1}{|b||a+b|}+\frac{1}{|a||b|}
\end{equation}
moreover
$
|m| \le  \int_s^t \left|\int_s^\sigma  e^{i b \sigma_1}
  d\sigma_1\right| d\sigma \le |b|^{-1}|t-s|  
$
and symmetrically $|m| \le |t-s| |a|^{-1}$. These last bounds imply
that  $|m| \le |t-s||ab|^{-1/2}$ and, since $2\gamma \in [0,1]$,  interpolating between this bound and
eq.~(\ref{eq:m-bound-1}) we get
\begin{equation}
  \label{eq:m-bound-3}
 |m| \le \frac{|t-s|^{2\gamma}}{|ab|^{\gamma}}
 \left(\frac{1}{|b||a+b|}+\frac{1}{|a||b|}\right)^{1-2\gamma} \le 
\frac{|t-s|^{2\gamma}}{|b|^{\gamma}|a|^{1-\gamma}|a+b|^{1-2\gamma}}+\frac{|t-s|^{2\gamma}}{|a|^{1-\gamma}|b|^{1-\gamma}}
\end{equation}
When $a+b=0$ then
$
m =   (1-e^{ia(t-s)})/a^2
$
which gives 
$
|m|\lesssim \frac{|t-s|^{2\gamma}}{|a|^{2-2\gamma}} \lesssim \frac{|t-s|^{2\gamma}}{|a|^{1-\gamma}|b|^{1-\gamma}} .
$
\end{proof}

\subsection{Proof of Lemma~\ref{lemma:X}}

\begin{proof} 
Estimating the norm of $X_{ts}$ is equivalent to estimating in $\FF L^p$ the norm of the operator $\YY_2(\Psi^1_{ts})\in \LL^2\FF L^p$ where $\Psi^1_{ts}: \RR^3\to \bbC$ is the multiplier
$$
\Psi^1_{ts}(k,k_1,k_2) =  \frac{|k|^{\alpha}}{|k_1|^\alpha |k_2|^\alpha} \frac{e^{-i 3 k k_1 k_2
      s}-e^{-i 3 k k_1 k_2 t}}{k_1 k_2} I_{kk_1k_2 \neq 0 }
$$ 
We must prove that, for $(\gamma,\alpha) \in \DD$ we have 
$|\YY_2(\Psi^1_{ts})|_{\LL^2 \FF L^p} \lesssim |t-s|^\gamma
$.
To begin we bound 
$$
|\Psi^1_{ts}(k,k_1,k_2)| \lesssim |t-s|^\gamma |k|^{2\gamma-1} \left(\frac{|k|}{|k_1| |k_2|}\right)^{1+\alpha-\gamma}  I_{kk_1k_2 \neq 0 }
$$
for any $\gamma \in[0,1]$.
We first prove a bound for this operator for $p=1,2,\infty$ and then conclude by interpolation.
For $p=1$, using Lemma~\ref{lemma:op-bound},  we have to bound 
$$
|\YY_2(\Psi^1_{ts})|_{\LL^2 \FF L^1}\lesssim
|t-s|^\gamma \sup_{k,k_1} 
 |k|^{2\gamma-1} \left(\frac{|k|}{|k_1| |k_2|}\right)^{1+\alpha-\gamma} I_{kk_1k_2 \neq 0 } $$
which is finite if $2\gamma-1 \le 0$ and $1+\alpha-\gamma\ge 0$ since from $k_1+k_2 = k$ we infer that $\max(|k_1|,|k_2|) \ge |k|/3$.
For $p=\infty$ using again Lemma~\ref{lemma:op-bound} we get
$$
|\YY_2(\Psi^1_{ts})|_{\LL^2 \FF L^\infty}\lesssim
|t-s|^\gamma \sup_{k \neq 0} |k|^{\alpha+\gamma} \sum_{k_1 \neq 0,k} 
  \frac{1}{|k_1|^{1+\alpha-\gamma} |k_2|^{1+\alpha-\gamma}}
 $$
which is finite if $\gamma\le 1/2$ and $1+\alpha-\gamma>1$ or by Cauchy-Schwartz if $\alpha+\gamma\le 0$ and $1/2+\alpha-\gamma > 0$. 
Interpolating the $p=1$ and $p=\infty$ bounds we get that for $p=2$ the norm is finite also for $\gamma=1/2$ and $\alpha > 0$. But when $p=2$ we have at our disposal another inequality, namely,
$$
|\YY_2(\Psi^1_{ts})|_{\LL^2 \FF L^2}^2 \lesssim
|t-s|^{2\gamma} \sum_{k\neq 0} |k|^{2(2\gamma-1)}\sup_{k_1 \neq 0,k} \left(\frac{|k|}{|k_1| |k_2|}\right)^{2(1+\alpha-\gamma)} 
$$
if $1+\alpha-\gamma\ge 0$ and $\gamma < 1/4$ this quantity is finite since can be bounded by
$
 \lesssim
|t-s|^{2\gamma} \sum_k |k|^{2(2\gamma-1)} \lesssim 1  .
$
Putting all together we get that for $p=2$ the norm is bounded when $\gamma < 1/4$ and $1+\alpha-\gamma\ge 0$ or when $\gamma\in[1/4,1/2]$ and $\alpha \ge -3/2+3 \gamma$.
Again by interpolation we obtain the bound for all $1 \le p \le +\infty$ and $(\gamma,\alpha) \in \DD$.
\end{proof}

\subsection{Proof of Lemma~\ref{lemma:X2}}

\begin{proof}
We start the argument  as in the proof of Lemma~\ref{lemma:X}. The Fourier transform of $X^2$ reads
\begin{equation*}
  \begin{split}
\FF X^2_{ts}&(\varphi_1,\varphi_2,\varphi_3)(k) 
 =
- 2 \sum'_{k_1,k_{21}}  \frac{k k_2}4  \int_s^t  \int_s^\sigma e^{-i 3
  k k_1 k_2 \sigma
- i 3 k_2 k_{21} k_{22} \sigma_1}
 \, d\sigma   d\sigma_1\, \wh \varphi_1(k_1)  \wh \varphi_2(k_{21})  \wh \varphi_3(k_{22})  
 \end{split}
\end{equation*}
then letting
$$
\Phi_{ts}^2(k,k_1,k_{21},k_{22}) = -   \frac{k k_2}2 I_{k k_1 k_2  k_{21} k_{22}\neq 0} \int_s^t  \int_s^\sigma e^{-i 3
  k k_1 k_2 \sigma
- i 3 k_2 k_{21} k_{22} \sigma_1} 
 \, d\sigma   d\sigma_1  ;
$$  
we have $X_{ts}^2 =\YY_3(\Phi_{ts}^2)$. Indeed
 it is enough to restrict the sums over the set of $k_i$ for  which $k k_1 k_2  k_{21} k_{22}\neq 0$.
The multiplier $\Phi^2_{ts}$
has two different behaviors depending on the quantity $h = k k_1 k_2 + k_2 k_{21} k_{22}$ being zero or not. We let $\Phi^2_{ts} = \hat\Phi^2_{ts}+  \breve \Phi^2_{ts}$ where 
$$
\breve \Phi^2_{ts} = \frac{k k_2}{i 3 k_2 k_{21} k_{22}} I_{h=0,k k_1 k_2  k_{21} k_{22}\neq 0} (t-s) = - \frac{t-s}{i 3 k_1} I_{h=0,k k_1 k_2  k_{21} k_{22}\neq 0} 
$$
$$
 = i \frac{t-s}{3 k_1} (I_{k=k_{21},k_1=-k_{22},k_2 \neq 0}+I_{k=k_{22},k_1=-k_{21},k_2 \neq 0})
$$
since the factorization $h=k_2(k_1+k_{21})(k-k_{21})$ implies that in the expression for $\breve \Phi^2_{ts}$ the only relevant contributions come from the case where $k=k_{21}$ and $k_1=-k_{22}$ or $k=k_{22}$ and $k_1 =-k_{21}$ and thus defining the
operators $\hat X^2_{ts}=\YY_3(\hat\Phi^2_{ts})$ and $\breve X^2_{ts}=\YY_3(\breve\Phi^2_{ts})$ we have $X^2_{ts} = \hat X^2_{ts}+\breve X^2_{ts}$ with  $\der \breve X^{2}_{tus} = 0$.

\paragraph{Bound for $\breve X^2$.}
When $1 \le p \le 2$ the operator $\breve X^2_{ts}$ is bounded in $\LL^3 \FF L^{\alpha,p}$ when $\alpha \ge -1/2$.
While when $p> 2$ is bounded  if 
$
[\sum_{k\neq 0} (|k|^{-1-2\alpha})^{p/(p-2)}]^{(p-2)/p}
$ is finite. This happens if $\alpha >-1/p$. In both cases we have $|\breve X^2_{ts}|_{\LL^3 \FF L^{\alpha,p}} \lesssim |t-s|^2$. Note that, due to some cancellations, the symmetric part of $\breve X^2$ is more regular and is bounded in any $\FF L^{\alpha,p}$ as soon as $\alpha \ge -1/2$. Indeed
\begin{equation}
\label{eq:cancellation-breve-X2}
\FF
\breve X^{2}_{ts}(\varphi, \varphi, \varphi)(k) = 2 (t-s) \sum_{k_1} \frac{ I_{0<|k|,|k_1|,|k_2|}}{3 i k_1} \wh \varphi (k_1) \wh \varphi (-k_1) \wh \varphi (k)
= -  \frac{2 (t-s) }{3 i k} \wh \varphi (k) \wh \varphi (-k) \wh \varphi (k).
\end{equation}
 We will not exploit here this better regularity since $\hat X^2$ will in any case limit the overall regularity of $X^2$ to $\alpha > -1/p$ when $3\gamma > 1$.

\paragraph{Bound for $\hat X^2$.}
For any $0\le \gamma \le 1/2$, Lemma~\ref{lemma:M-bounds} gives
$$
|\hat \Phi^2_{ts}(k,k_1,k_{21},k_{22})|\lesssim (A_1+A_2)  |t-s|^{2\gamma} I_{k_1,k_{21},k_{22} \neq 0,k}
$$
where 
$
A_1 = |k| |k_2| |k k_1 k_2|^{\gamma-1} |k_2 k_{21} k_{22}|^{\gamma-1}
$
and
$$
A_2 = |k|^{1-\gamma}  |k_1|^{-\gamma}|k_{22}|^{\gamma-1}|k_{21}|^{\gamma-1}
|k-k_{21}|^{2\gamma-1} 
|k-k_{22}|^{2\gamma-1} |k-k_1|^{2\gamma-1} .
$$
So in order to bound $\hat X^2_{ts}$ in $\FF L^{\alpha,p}$ with a quantity of order $|t-s|^{2\gamma}$ it will be enough to bound separately
the multipliers
$$
\Theta_1 =  A_1 |k|^{\alpha} |k_1|^{-\alpha}
|k_{21}|^{-\alpha} |k_{22}|^{-\alpha} I_{k_1,k_{21},k_{22} \neq 0,k}
$$
and
$$
\Theta_2 =  A_2 |k|^{\alpha} |k_1|^{-\alpha}
|k_{21}|^{-\alpha} |k_{22}|^{-\alpha}  I_{k_1,k_{21},k_{22} \neq 0,k} 
$$
$$
= I_{k_1,k_{21},k_{22} \neq 0,k}\left( \frac{|k|}{|k_1||k_{21}||k_{22}|}\right)^{1+\alpha-\gamma}
\left( \frac{|k_1|}{|k-k_1||k-k_{21}||k-k_{22}|}\right)^{1-2\gamma}
$$
in $\FF L^p$.
The expression $A_1$ is nicely factorized in its dependence on the
couple $k_{21},k_{22}$ and  the multiplier $\Theta_1$ can be easily bounded
using (twice and iteratively) the same arguments as in
Lemma~\ref{lemma:X}.
We will concentrate on the multiplier $\Theta_2$ which requires  different
estimates. 

 Since $k_1+k_{21}+k_{22}=k$ and $(k-k_1)-(k-k_{21})-(k-k_{22}) = -2 k_1 $ at least one of $|k_1|,|k_{21}|,|k_{22}|$ is larger than $|k|/4$ and similarly one of $|k-k_1|,|k-k_{21}|,|k-k_{22}|$ is larger than $|k_1|/2$. Using this fact,
when $1+\alpha-\gamma \ge 0$ and $1-2\gamma \ge0$,  we have
$$
\Theta_2 \lesssim I_{k_a,k_{b},k_{c} \neq 0,k}
\left( \frac{1}{|k_a||k_b|}\right)^{1+\alpha-\gamma}
\left( \frac{1}{|k-k_d||k-k_e|}\right)^{1-2\gamma}
$$
where $(a,b,c)$ and $(d,e,f)$ are two permutations of $(1,21,22)$ (depending on $k_1,k_{21},k_{22}$) so that we must have $\{a,b\} \cap \{d,e\} \neq \emptyset$.
It is clear then that $\Theta_2 \lesssim 1$ so that $\YY_3(\Theta_2)$ is bounded in $\LL^3\FF L^{\alpha,1}$ when
$1+\alpha-\gamma \ge 0$ and $1-2\gamma \ge0$.
The bound in $\LL^3 \FF L^{\alpha,\infty}$ will follow by showing that
$
\sup_k \sum_{k_1,k_{21}} \Theta_2(k,k_1,k_{21},k_{22})
$
is finite. If $1+\alpha-\gamma > 1$ and $1-2\gamma \ge 0$ we can bound this quantity by
$$
\sup_k \sum_{k_1,k_2 \neq 0,k} \left( \frac{1}{|k_1||k_{2}|}\right)^{1+\alpha-\gamma} \left( \frac{1}{|k-k_1|}\right)^{1-2\gamma} \lesssim 1 .
$$
When $1+\alpha-\gamma \ge 0$ a bound in any $\LL^3 \FF L^{\alpha,p}$ for $p\in [1,+\infty]$ is obtained from a bound of
$$
\sup_k \sum_{k_1,k_{21}} 
[ \Theta_2(k,k_1,k_{21},k_{22}) ]^{q} \lesssim \sum_{n_1,n_2}  \left( \frac{1}{|n_1|}\right)^{q(1+\alpha-\gamma)}
\left( \frac{1}{|n_1||n_2|}\right)^{q(1-2\gamma)} .
$$
 This last 
quantity is then finite for  any $\gamma < 1/(2p)$.
By interpolation of these various estimates we obtain again that $\YY_3(\Theta_2)$ is bounded in $\LL^3 \FF L^{\alpha,p}$ for any $(\gamma,\alpha) \in \DD$. Note that the extra factors in the majorizing sums can be used to show that actually in the interior of the region $\DD$ the operator $X^2$ in bounded from $\FF L^{\alpha,p}$ to $\FF L^{\alpha+\eps,p}$ for some $\eps>0$. 
\end{proof}

\subsection{Proof of Lemma~\ref{lemma:gal}}

\begin{proof}
Let us work out the details for $X$, first.
The Fourier multiplier associated to $X^{(N)}_{ts}$ differs from that of $X_{ts}$ for a factor of the form
$
I_{|k|\le N,|k_1|\le N,|k_2|\le N}
$ so the difference $X_{ts}-X^{(N)}_{ts}$ has a Fourier multiplier containing the factor
$$
I_{kk_1k_2 \neq 0}(1-I_{|k|\le N,|k_1|\le N,|k_2|\le N}) \le I_{kk_1k_2 \neq 0}( I_{|k_1| > N}+I_{|k_2| > N}+I_{|k| > N} )\lesssim \frac{|k k_1 k_2|^\eps}{N^\eps} I_{kk_1k_2 \neq 0}
$$
for any $\eps > 0$.
To prove the convergence in $\CC^\gamma_2 \FF L^{\alpha,p}$ is then enough to prove that $X$ is bounded in the space $\CC_2^\gamma \LL((\FF L^{\alpha-\eps,p})^2,\FF L^{\alpha+\eps,p})$ for some $\eps > 0$. When $(\gamma,\alpha)$ is in the interior of $\DD$, exploiting the gaps in the inequalities involved in the proof of Lemma~\ref{lemma:X}, it is not difficult to show that this is the case for sufficiently small $\eps > 0$.
Then
\begin{equation*}
  \begin{split}
 | (X_{ts}-X^{(N)}_{ts})(\varphi_1,\varphi_2)|_{{\FF L^{\alpha,p}}} & \lesssim N^{-\eps} |\varphi_1|_{{\FF L^{\alpha,p}}} |\varphi_2|_{{\FF L^{\alpha,p}}} .
  \end{split}
\end{equation*}
A similar argument works for the convergence of $\hat X^{(N),2}$ to $\hat X^{(2)}$  since in this case the multiplier associated to the difference $\hat X^2-\hat X^{(N),2}$ has a factor of the form
$$
I_{kk_1k_2 k_{21} k_{22} \neq 0}(1-I_{|k|,|k_1|,|k_2|,|k_{21}|,|k_{22}|\le N}) 
\lesssim \frac{|k k_1 k_2 k_{21} k_{22}|^\eps}{N^\eps} I_{kk_1k_2 k_{21} k_{22} \neq 0} 
\lesssim \frac{|k k_1 k_{21} k_{22}|^{2\eps}}{N^\eps} I_{kk_1k_2 k_{21} k_{22} \neq 0} 
$$
where we used the fact that $k_2 = k_{21}+k_{22}$ to remove it from the r.h.s. Again an inspection of the proof of the regularity of $\hat X^2$ confirms that in the interior of $\DD$ we have a small gain of regularity which can be converted in an estimate of $\hat X^2$ in $\CC^{2\gamma} \LL^3((\FF L^{\alpha-\eps,p})^2;\FF L^{\alpha+\eps})$ and obtain the norm convergence of $\hat X^{(N),2}$ to $\hat X^2$ with a polynomial speed:
$
 |\hat X^2_{ts}-\hat X^{(N),2}_{ts}|_{{\LL^3 \FF L^{\alpha,p}}} \lesssim N^{-\eps/2} |t-s|^{2\gamma}
$
for some $\eps > 0$.
\end{proof}

\section{Proof of Theorem~\ref{thm:main}}
\label{proof:thm:main}
\begin{proof} Fix $\eta \in (0,\gamma)$ and
 define
a map $\Gamma: \QQ_\eta \to \QQ_\eta$ by $z = \Gamma(y)$ where
\begin{equation}
  \label{eq:gamma}
z_0 = y_0, \qquad \der z = X(y,y) + X^2(y,y,y) + z^\flat
\end{equation}
with
$$
z^\flat = \Lambda[ 2 X(y^\sharp,y) 
+ X(\der y,\der y) 
 + X^2(\der y,y,y) +
X^2(y,\der y,y) + X^2( y,y,\der y)   ].
$$
and we set $z' = y$ and $z^\sharp = X^2(y,y,y) + z^\flat$. For the map
$\Gamma$ to be well defined it is enough that $(\gamma,\alpha) \in
\DD$ and that $\eta > 1/3$. Indeed if $y \in \QQ_\eta$ then all the
terms in the argument of the $\Lambda$-map belongs to
$\CC_2^{3\eta} \FF L^{\alpha,p}$.
 
We will proceed in three steps. First we will prove that there exists
$T_*\in(0,1]$  depending only on $|v_0|$ such that for any $T\le T_* $, $\Gamma$ maps a closed ball $B_T$ of $\QQ_{\eta}$ in
itself. Then we will prove that if $T$ is sufficiently small, $\Gamma$
is actually a contraction on this ball proving the existence  of a
fixed-point. Finally uniqueness will follow from a standard argument.

\paragraph{Boundedness.}
The key observation to prove (below) that $\Gamma$
is a contraction for $T$ sufficiently small is that $\Gamma(y)$
belongs actually to $\QQ_\gamma$, indeed $z' = y \in \CC_1^\gamma
 \FF L^{\alpha,p}$ and $z^\sharp \in \CC_2^{2\gamma}  \FF L^{\alpha,p}$. 

Let us assume that $T \le 1$.
First we will prove that
\begin{equation}
  \label{eq:bound-on-Gamma}
  d_{\QQ,\gamma}(\Gamma(y),0) = d_{\QQ,\gamma}(z,0) \le C (1+d_{\QQ,\eta}(y,0))^3.
\end{equation}
Since $3\eta > 1$ using Prop.~\ref{prop:Lambda} we have
\begin{equation}
\label{eq:flat-bound}
  \begin{split}
\|z^{\flat}\|_{3\eta} & \le C \|2 X(y^\sharp,y) 
+ X(\der y,\der y) 
 + X^2(\der y,y,y) +
X^2(y,\der y,y) + X^2( y,y,\der y) \|_{3\eta}
\\ & \le C \left[ 2\| X(y^\sharp,y) \|_{3\eta}
+ \|X(\der y,\der y)\|_{3\eta} 
\right. \\ & \qquad \left. + \|X^2(\der y,y,y)\|_{3\eta} +
\|X^2(y,\der y,y)\|_{3\eta} + \|X^2( y,y,\der y) \|_{3\eta}\right]
  \end{split}
\end{equation}
The first term is bounded by
$
\| X(y^\sharp,y) \|_{3\eta} \le \|X\|_{\eta} \|y^\sharp\|_{2\eta} \|y\|_{\infty}
$
where the H\"older norm $\|X\|_{\eta}$ is considered with respect to
the operator norm on $\LL(( \FF L^{\alpha,p})^2, \FF L^{\alpha,p})$. Now
$\|X\|_\eta \le \|X\|_\gamma$ since we assume $\eta < \gamma$, while
$\|y^\sharp\|_{2\eta} \le d_{\QQ,\eta}(y,0)$ and finally 
remark that $\|y\|_{\infty} = |y_0|+ \|y\|_{\eta} \le
d_{\QQ,\eta}(y,0)$ since we assumed $T \le 1$. So  
$
\| X(y^\sharp,y) \|_{3\eta} \le C d_{\QQ,\eta}(y,0)^2
$.
 Each of the
other terms in eq.~(\ref{eq:flat-bound}) can be bounded similarly.
 Then 
\begin{equation}
\label{eq:flat-bound-final}
  \begin{split}
\|z^{\flat}\|_{3\eta} & \le C (1+d_{\QQ,\eta}(y,0))^3
  \end{split}
\end{equation}
This implies 
\begin{equation*}
  \begin{split}
\|z^\sharp\|_{2\gamma} & = \|X^2(y,y,y) + z^\flat\|_{2\gamma}   
 \le \|X^2(y,y,y)\|_{2\gamma} + \|z^\flat\|_{2\gamma}
\\ & \le C \|y\|_{\infty}^3 + \|z^\flat\|_{2\gamma}     
  \end{split}
\end{equation*}
By hypothesis we have $2\gamma < 1 < 3\eta$ so $\|z^\flat\|_{2\gamma}
\le  \|z^\flat\|_{3\eta}$ and again we obtain
$
\|z^\sharp\|_{2\gamma} \le C (1+d_{\QQ,\eta}(y,0))^3
$.
Finally we have to bound $z'$ and $z$:
\begin{equation*}
  \begin{split}
\|z'\|_{\gamma} & = \|y\|_{\gamma} \le \|X(y',y')+y^\sharp\|_{\gamma} \le
C \|y'\|_{\infty}^2 + \|y^\sharp\|_{2\gamma}
\\  & \le
C (|y_0|^2 + \|y'\|_{\gamma}^2) + \|y^\sharp\|_{2\gamma}
\le C (1+d_{\QQ,\eta}(y,0))^2    
  \end{split}
\end{equation*}
Similarly we can bound $z$ as $\|z\|_\gamma \le C
(1+d_{\QQ,\eta}(y,0))^3$.
Putting these bounds together we prove eq.~(\ref{eq:bound-on-Gamma}).
Let $\eps = \gamma-\eta > 0$ and observe that
\begin{equation*}
  \begin{split}
    d_{\QQ,\eta}(\Gamma(y),0) &= |z_0| +
    \|z\|_\eta+\|z'\|_{\eta}+\|z^\sharp\|_{2\eta} 
\\ & \le |z_0| +
     T^\eps\left[ \|z\|_\gamma+\|z'\|_{\gamma}+\|z^\sharp\|_{2\gamma}\right]
\\ & \le |y_0| +
     T^\eps d_{\QQ,\gamma}(\Gamma(y),0)
  \end{split}
\end{equation*}
So 
\begin{equation}
\label{eq:boundedness}
    d_{\QQ,\eta}(\Gamma(y),0) \le |y_0|+ C T^\eps (1+d_{\QQ,\eta}(y,0))^3
\end{equation}
Let $B_T = \{ y \in \QQ_{\eta} : d_{\QQ,\eta}(y,0) \le b_T, y_0 = v_0
\}$ be a closed ball of elements in $\QQ_{\eta}$ with initial
condition $v_0$. Let $b_T > 0$ be the solution of the algebraic equation
$
b_T = |v_0|+ C T^\eps (1+b_T)^3
$
which exists whenever $T \le T_*(|v_0|)$ where $T_*(|v_0|)$ is a
strictly positive time which depends only on $|v_0|$. 
If $T \le T_*(|v_0|)$ and $y \in B_T$ then by eq.~(\ref{eq:boundedness}) we have
$
d_{\QQ,\eta}(\Gamma(y),0) \le |v_0| + C T^\eps (1+b_T)^3 = b_T
$    
which implies that $\Gamma$ maps $B_T$ onto itself.

\paragraph{Contraction.}
By arguments similar to those used above we can prove that if $y^{(1)},
y^{(2)}$ are two elements of $\QQ_{\eta}$ then
\begin{equation}
  \label{eq:Gamma-contraction}
  d_{\QQ,\gamma}(\Gamma(y^{(1)}),\Gamma(y^{(2)})) \le C d_{\QQ,\eta}(y^{(1)},y^{(2)}) (1+d_{\QQ,\eta}(y^{(1)},0)+d_{\QQ,\eta}(y^{(2)},0))^2 
\end{equation}
Then 
\begin{equation}
  \label{eq:Gamma-contraction-2}
   d_{\QQ,\eta}(\Gamma(y^{(1)}),\Gamma(y^{(2)})) \le
  |y^{(1)}_0-y^{(2)}_0| 
  + C T^\eps  d_{\QQ,\eta}(y^{(1)},y^{(2)}) (1+ d_{\QQ,\eta}(y^{(1)},0)+ d_{\QQ,\eta}(y^{(2)},0))^2 
\end{equation}
Taking $y_1,y_2 \in B_T$ we have
\begin{equation}
  \label{eq:Gamma-contraction-3}
  d_{\QQ,\eta}(\Gamma(y^{(1)}),\Gamma(y^{(2)})) \le
   C (1+2 b_T)^2 T^\eps d_{\QQ,\eta}(y^{(1)},y^{(2)})
\end{equation}
Let $T_{**} < T_*$ and such that  $C (1+2 b_T)^2 T_{**}^\eps < 1$, eq.~(\ref{eq:Gamma-contraction-3}) implies that
$\Gamma$ is a strict contraction on $B_{T_{**}}$ with a unique fixed-point
(in $B_{T_{**}}$).
A standard argument allows to extend the solution 
to the larger time interval $T_*$.

\paragraph{Uniqueness.} Now assume to have two solutions 
$v^{(1)},v^{(2)} \in \QQ_\eta$ of eq.~(\ref{eq:abs-step-2}), i.e. such
that $v^{(i)} = \Gamma(v^{(i)})$, $v^{(i)}_0 = v_0$ for
$i=1,2$. Eq.~(\ref{eq:Gamma-contraction}) applied to the distance between $\Gamma(v^{(1)})
$ and $\Gamma(v^{(2)})$ implies
\begin{equation*}
 d_{\QQ,\gamma}(v^{(1)},v^{(2)}) \le C  d_{\QQ,\eta}(v^{(1)},v^{(2)}) (1+ d_{\QQ,\eta}(v^{(1)},0)+ d_{\QQ,\eta}(v^{(2)},0))^2   
\end{equation*}
If we denote $\QQ_{\eta,T}$ the space $\QQ_{\eta}$ considered for
functions on $[0,T]$ and with $ d_{\QQ,\eta,T}(\cdot,\cdot)$ the
corresponding distance, then we have for $S \le T$
$
 d_{\QQ,\eta,S}(y,0) \le   d_{\QQ,\eta,T}(y,0)
$ for any $y \in \QQ_{\eta,T}$. Then
\begin{equation*}
 d_{\QQ,\eta,S}(v^{(1)},v^{(2)}) \le C  d_{\QQ,\eta,S}(v^{(1)},v^{(2)}) (1+ d_{\QQ,\eta,S}(v^{(1)},0)+ d_{\QQ,\eta,S}(v^{(2)},0))^2   
\end{equation*}
and so
\begin{equation*}
 d_{\QQ,\eta,S}(v^{(1)},v^{(2)}) \le C S^\eps  d_{\QQ,\eta,S}(v^{(1)},v^{(2)}) (1+ d_{\QQ,\eta,S}(v^{(1)},0)+ d_{\QQ,\eta,S}(v^{(2)},0))^2   
\end{equation*}
so when $S$ is small enough that $C S^\eps
(1+ d_{\QQ,\eta,S}(v^{(1)},0)+ d_{\QQ,\eta,S}(v^{(2)},0))^2  \le 1/2  $ we
must have that $ d_{\QQ,\eta,S}(v^{(1)},v^{(2)}) = 0$ and then
$v^{(1)}_t = v^{(2)}_t$ for any $t \in [0,S]$. By repeating this
argument  we conclude that the two solution coincide in the whole interval
$[0,T]$.

\paragraph{Lipschitz continuity.}
Modulo some technicalities due to the local character of the solution, the proof of the Lipschitz continuity of the map $\Theta$ follows the line of the finite dimensional situation. See~\cite{MR2091358} for details.
\end{proof}


\end{document}